\documentclass[smallextended,referee,envcountsect]{svjour3}
\smartqed
\usepackage{amsmath,amsfonts,amssymb,amscd}
\usepackage{bbm}
\usepackage{enumerate}
\usepackage{enumitem}
\usepackage{url}
\usepackage{cite}
\usepackage{fullpage}
\usepackage{fancyhdr}
\usepackage{color}
\usepackage{float}
\usepackage{soul}
\usepackage{graphicx}

\usepackage{hyperref}

\usepackage[margin=1in]{geometry}

\renewcommand\theenumi{(\roman{enumi})}
\renewcommand\theenumii{(\alph{enumii})}

\newcommand{\inner}[2]{\langle{#1},{#2}\rangle}

\usepackage{cite}

\newcommand{\argmin}{{\rm arg}\!\min}

\newcommand{\la}{\langle}
\newcommand{\ra}{\rangle}
\newcommand{\ve}{\varepsilon}

\newcommand{\dom}{{\bf dom\,}}

\newcommand{\prox}{{\bf prox}}
\newcommand{\nexto}{\kern -0.54em}

\newcommand{\dZ}{{\cal Z \kern -0.7em Z}}
\newcommand{\dC}{{\rm\hbox{C \kern-0.8em\raise0.2ex\hbox{\vrule height5.4pt width0.7pt}}}}
\newcommand{\dQ}{{\rm\hbox{Q \kern-0.85em\raise0.25ex\hbox{\vrule height5.4pt width0.7pt}}}}

\usepackage{epsfig}
\usepackage{latexsym}

\def\ve{\varepsilon}

\newcommand{\bi}{\begin{itemize}}
\newcommand{\ei}{\end{itemize}}
\newcommand{\ba}{\begin{array}}
\newcommand{\ea}{\end{array}}

\journalname{JOTA}

\begin{document}

\title{A Relative Inexact Proximal Gradient Method with an Explicit  Linesearch }


\author{Yunier\ Bello-Cruz \and Max L.N. Gon\c calves\and Jefferson G. Melo\and Cassandra Mohr}

\institute{Yunier\ Bello-Cruz  \and Cassandra Mohr \at
Department of Mathematical Sciences, Northern Illinois University. Watson Hall 330, DeKalb, IL, USA - 60115.\\
              yunierbello@niu.edu and  cmohr3@niu.edu.
           \and
             Max L.N. Gon\c calves \and   Jefferson G. Melo  \at
     Institute of Mathematics and Statistics, Federal University of Goias, Campus II- Caixa
    Postal 131, CEP 74001-970, Goi\^ania-GO, Brazil.\\
             maxlng@ufg.br and jefferson@ufg.br
}

\date{Received: date / Accepted: date\\[2mm]Communicated by Nobuo Yamashita}

\maketitle

\begin{abstract}
This paper presents and investigates an inexact proximal gradient method for solving composite convex optimization problems characterized by an objective function composed of a sum of a full-domain differentiable convex function and a non-differentiable convex function. 
We introduce an explicit line search applied specifically  to the differentiable component of the objective function, requiring only a relative inexact solution of the proximal subproblem per iteration.
We prove the convergence of the sequence generated by our scheme and establish its iteration complexity, considering both the functional values and a residual  associated with first-order stationary solutions. Additionally, we provide numerical experiments to illustrate  the practical efficacy of our method.
\end{abstract}
\keywords{ Linesearch \and  Iteration
complexity \and  Nonsmooth and convex optimization problem\and  Proximal
gradient method \and  Relative error rule}
\subclass{65K05 \and 
90C25 \and 90C30}


\section{Introduction}

In this paper, our focus is on addressing nonsmooth convex optimization problems characterized by the following formulation: 
\begin{equation}\label{main-prob}
\min \, F(x) := f(x) + g(x) \quad \text{subject to} \quad x \in \mathbbm{E},
\end{equation}
where \( \mathbbm{E} \) represents a finite-dimensional Euclidean space. The function \( g : \mathbbm{E} \to \overline{\mathbbm{R}} := \mathbbm{R} \cup \{+\infty\} \) is a nonsmooth, proper, and lower semicontinuous convex function. The function \( f : \mathbbm{E} \to \mathbbm{R} \) is a continuously differentiable and convex function. Throughout the paper, we denote the optimal value and solution set of problem~\eqref{main-prob} by $F_*$ and $S_*$, respectively. From now on, we assume that $S_*\neq \emptyset$.

Proximal gradient methods effectively solve optimization problems such as \eqref{main-prob}. The main step in the proximal gradient method involves evaluating the proximal operator $\prox_g:\mathbbm{E} \to \dom g:=\{x\in\mathbbm{E} \mid g(x)<+\infty\}$ defined as follows:
\begin{equation} \label{section1_1.1}
\prox_{g}(x) := \underset{y \in \mathbbm{E}}{\arg\min} \; \left\{ g(y) + \frac{1}{2} \left \| y - x \right \|^2 \right\},
\end{equation} where the norm, $\|\cdot\|$, is induced by the inner product of $\mathbbm{E}$, $\la \cdot, \cdot\ra$, as $\|\cdot\|:=\sqrt{\la \cdot, \cdot\ra}$. The proximal operator is well-known as a full-domain, firmly nonexpansive operator. These useful properties, together with the descent property of the gradient step, establish the foundation for the convergence and complexity analysis for proximal gradient iterations.

\begin{center}\fbox{\begin{minipage}[b]{0.97\textwidth}
\noindent {\bf Proximal Gradient Method (PGM)}

Let $x_0\in \dom g$. Compute $\lambda_{k}>0$ and 
	\begin{equation}\label{bar-xk-PGM}
	\tilde x_{k} := \prox_{\lambda_{k} g}\left(x_{k}-\lambda_{k}\nabla f(x_{k})\right).
	\end{equation}
	Choose $\beta_{k}\in (0,1]$ and compute
	 \begin{equation}\label{xk-PGM}
	x_{k+1} := x_{k}+\beta_{k} (\tilde x_{k}-x_{k}).
	\end{equation}
\end{minipage}}\end{center}
The coefficients $\lambda_{k}$ and $\beta_{k}$, referred to as stepsizes, can be determined based on backtracking linesearch procedures. Such strategies are essential whenever the global $L$-Lipschitz continuity of the gradient of $f$ fails or even when computing an acceptable upper bound for $L$ is challenging. This situation is often encountered in numerous applications; for instance, in inverse problems based on non-Euclidean norms \cite{Bredies:2009,Schuster:2012} or Bregman distances such as the Kullback-Leibler divergence \cite{Salzo:2014,Bonettini:2012,Csiszar:1991,Vardi:1985}. Moreover, even when $L$ is known, linesearches may allow for longer steps toward the solution by using local information at every iteration.

There are several possible choices for these stepsizes, each impacting the algorithm's performance in different ways; see, for instance, \cite{Salzo:2017a, BelloCruz:2016a,BelloCruz:2016c,Bello-Cruz:2021}. It is important to note that in order to compute the stepsize $\lambda_{k}$ using a backtracking linesearch at each iteration $k$, the proximal operator may need to be evaluated multiple times within the procedure. Conversely, the stepsize $\beta_{k}$ can be selected by evaluating the proximal operator only once per iteration. In this context, we will refer to explicit linesearch to describe a backtracking procedure that determines $\beta_{k}$ after setting $\lambda_{k}$ as a constant for all $k$.
This type of explicit strategy, first presented in \cite{BelloCruz:2016a}, is particularly advantageous, especially in cases where evaluating the proximal operator is challenging.  
The function $g$ is often complex enough that the corresponding proximal operator lacks an analytical solution. In such cases, an ad-hoc algorithm should be employed to evaluate the proximal operator inexactly. For instance, consider $\mathbb{E} = \mathbbm{R}^n$, and let $g: \mathbbm{R}^n \to \mathbbm{R}$ be defined  as $$g(x) = \| x \|_1 + \lambda\sum_{i<j} \max \{ | x_i | , | x_j | \},$$ with $\lambda > 0$, a form encountered in sparse regression problems with grouping, as discussed in \cite{Zhong:2012}. Similarly, in the context of matrix factorization, consider $\mathbb{E}=\mathbbm{R}^{n \times m}$ for the CUR-like factorization optimization problem \cite{Barre:2022,JMLR:v12:mairal11a}, where the goal is to approximate a matrix $W \in \mathbbm{R}^{m \times n}$ with $X \in \mathbbm{R}^{n \times m}$ having sparse rows and columns. In this case, $g:\mathbbm{R}^{n \times m}\to \mathbbm{R}$ is given  by $$g(X) = \lambda_{\text{row}}\sum_{i=1}^{n}\|X^{(i)}\|_2 + \lambda_{\text{col}}\sum_{j=1}^{m}\|X_{(j)}\|_2,$$ where $\lambda_{\text{row}}, \lambda_{\text{col}}>0$. This nondifferentiable term in problem~\eqref{main-prob} will be considered in the numerical illustrations of this paper.
 For further examples and discussions, see \cite{BelloCruz:2017a,Bello-Cruz:2023,Salzo:2012,Villa:2013,Schmidt:2011,Jiang:2012,Millan:2019,Barre:2022,BelloCruz:2022b}. Of course, there are some cases when the exact analytical solution of the proximal operator is available, such as when $g:\mathbbm{R}^n\to \mathbbm{R}$ as $g(x) = \| x \|_1$ or the indicator function of a simple convex and closed set, see, for instance, \cite{Beck:2009,Combettes:2005}.

Consequently, in practice, the evaluation of the proximal operator is often done inexactly. A basic inexactness criterion is to approximately evaluate the proximal operator using an exogenous sequence of error tolerances  which, in general, must be summable in order to guarantee the convergence of the proximal iterations. 
Such a diminishing sequence is a priori chosen without using any information that may be available along with the iterations. Usually, the restrictive summability condition forces the solution of the proximal subproblem to be increasingly accurate, often more than necessary; see, for instance, \cite{Jiang:2012,Salzo:2012,Aujol:2015,Villa:2013,Schmidt:2011}. In the past two decades, relative error criteria have been considered an effective  and practical way of controlling the inexactness in solving  proximal subproblems of several algorithms, including FISTA, ADMM, augmented Lagrangian, Douglas-Rachford, and {proximal gradient} methods;
refer to \cite{Eckstein:2013,Eckstein:2018,Alves:2020,Adona:2019,Adona:2020,Bonettini:2016,Bello-Cruz:2023,Millan:2019} for examples. Relative error criteria are often easy to verify in practice and have the advantage of exploiting the available information at a specific iteration, being, therefore,  an interesting alternative to the aforementioned exogenous sequences.

In this paper, we propose and analyze an inexact  proximal gradient method (PGM) for solving problem~\eqref{main-prob}.  We present a novel relative inexactness criterion for solving the proximal subproblems which somewhat resembles the ideas of relative error criteria introduced  in \cite{Monteiro:2010,Solodov:1999a}, but incorporates some new elements to control the objective function for the inexact solution. The proposed scheme requires only one inexact solution of the proximal subproblem per iteration, and the stepsizes are computed through a relaxed explicit linesearch procedure, applied specifically  to $f$, that takes  into account the residuals obtained from the proximal subproblem and  enables the iteration to address non-Lipschitz optimization problems effectively. We show that the sequence generated by our method converges to a solution of problem~\eqref{main-prob}. Moreover, we  establish its iteration complexity in terms of both the function values and the residuals associated with an approximate stationary solution.  We also present some numerical experiments to illustrate the performance of the proposed method. 

It is worth mentioning that the PGM proposed in this paper can be regarded as an inexact version of the ones analyzed in \cite{Salzo:2017a,BelloCruz:2016a}, where explicit linesearches were studied in the PGM setting.
On the other hand, an inexact version of the PGM with an explicit linesearch was also proposed in \cite{Bonettini:2016}. 
The authors developed a general framework based on variable metrics and established results on convergence rates for the sequence generated by their method.  
While, in the special case where the metric is fixed and chosen as the Euclidean one, their algorithm shares some similarities with ours, it employs  different type of linesearch and  inexact criteria for solving the proximal subproblem.
Further details on the differences between our algorithm and the method in \cite{Bonettini:2016} are discussed in Subsection~\ref{HPE-inexact-subsection}.


The paper is structured as follows: Section~\ref{prel} presents definitions, basic facts, and auxiliary results. 
The concept of an approximate solution for the proximal subproblem, along with the description of algorithms that can be employed to compute it, is detailed in Section~\ref{Inexact-Rules}. In this latter section, we discuss some works that analyze relative inexact proximal solution criteria related to ours.  
Section~\ref{search} introduces the inexact proximal gradient method with an explicit linesearch (IPG-ELS) and establishes some of its fundamental properties. Section~\ref{convSec} analyzes the full convergence of the sequence generated by the IPG-ELS method and establishes its iteration complexity bounds  in terms of functional values and a residual associated with the stationary condition for problem~\eqref{main-prob}.  
Some numerical experiments illustrating the behavior of the proposed scheme are
reported in Section~\ref{NumSec}. Finally, concluding remarks are provided in Section~\ref{concluding}.

\section{Preliminary }\label{prel}

In this section, we present some preliminary material and notations that will be used throughout this paper.

Let $g:\mathbbm{E}\to \overline{\mathbbm{R}}$ be a proper, lower semicontinuous, and convex function. For a given $\varepsilon \geq 0$, the $\varepsilon$-subdifferential of $g$ at $x \in \dom g=\{x\in\mathbbm{E} \mid g(x)<+\infty\}$, denoted by $\partial_\varepsilon g(x)$, is defined as 
\begin{equation}\label{def-epsSubdiff}
\partial_\varepsilon g(x) := \{v\in\mathbbm{E} \mid g(y) \ge g(x) + \langle v, y - x \rangle - \varepsilon, \ \forall y \in \mathbbm{E}\}.
\end{equation} When $x\notin \dom g$, we define $\partial_\varepsilon g(x)=\emptyset$. Any element $v\in \partial_\varepsilon g(x)$ is called an $\varepsilon$-subgradient. 
If $\varepsilon = 0$, then $\partial_0 g(x)$ becomes the classical subdifferential of $g$ at $x$, denoted by $\partial g(x)$.  It follows immediately from \eqref{def-epsSubdiff} that 
\begin{equation}\label{eps-monotonicity of subdiff g}
v \in \partial_\varepsilon g(y), \; u \in \partial g(x) \quad \text{implies} \quad \langle v - u, y - x \rangle \geq -\varepsilon.
\end{equation}
We present two useful auxiliary results for $\partial_\varepsilon g$.  The first one is the closedness of the graph of $\partial_\varepsilon g$ and the second is the so-called transportation formula; see Propositions 4.1.1 and   4.2.2  of \cite{Hiriart-Urruty:1993}.

\begin{proposition}[Closed Graph Property]\label{closed-graph} 
Let $(\varepsilon_k,x_k,v_k)_{k\in\mathbbm{N}}\subseteq \mathbbm{R}_+\times \mathbbm{E}\times \mathbbm{E}$ be a sequence converging to $(\varepsilon,x,v)$ with $v_k\in \partial_{\varepsilon_k} g(x_k)$ for all $k\in\mathbbm{N}$. Then, $v\in \partial_\varepsilon g(x)$.
\end{proposition}

\begin{proposition}[Transportation Formula]\label{transp}  With $x$ and $y$ in  $\dom g,$ let $v\in \partial g (y)$. Then $v\in \partial g_{\epsilon} (x)$ where $\epsilon= g(x)-g(y)-\inner{v}{x-y}\geq 0$.
\end{proposition}

We now introduce a concept of  approximate stationary solution to problem~\eqref{main-prob}, which can be seen as a specialization of the one presented in \cite[Eq. 1]{Monteiro:2010}. First, note that $\bar{x}$ is a solution to problem~\eqref{main-prob} if and only if $0 \in \nabla f(\bar{x}) + \partial g(\bar{x})$. The concept of approximate stationary solution relaxes this inclusion by introducing a residual $v$ and enlarging $\partial g$ using $\partial_\varepsilon g$.

\begin{definition}[$\eta$-Approximate Stationary Solution]\label{approxStationarySolution}
Given $\eta > 0$, an element $\tilde{x} \in \dom g$ is said to be an
$\eta$-approximate stationary solution to problem~\eqref{main-prob} with a residual pair $(v,\varepsilon)$ if 

\begin{equation}\label{def:eps-stationary-sol}
v \in \nabla f(\tilde{x}) + \partial_\varepsilon g(\tilde{x}), \quad \max\{\|v\|,\varepsilon\} \leq \eta.     
\end{equation}
\end{definition}

Next, we recall the definition of quasi-Fej\'er convergence, which is an important and well-known tool for establishing full convergence of gradient and proximal point type methods; see, for instance, \cite[Definition 1]{Iusem-Quasi-Fejer}.

\begin{definition}[Quasi-Fej\'er Convergence]\label{def-fejer}
Let $S$ be a nonempty subset of  $\mathbbm{E}$. A sequence $(x_{k})_{k\in\mathbbm{N}}\subseteq \mathbbm{E}$ is said to be quasi-Fej\'er convergent to $S$ if and only
if, for every $x \in S$, there exists a summable sequence $ (\delta_{k})_{k\in\mathbbm{N}} \subseteq \mathbbm{R}_+ $ such that
\begin{equation}\label{ineq:eps-Fejer}
\| x_{k+1}-x\|^2 \leq \| x_{k}-x\|^2 +  \delta_{k}, \qquad \forall k\in \mathbbm{N}.
\end{equation}
\end{definition}

The following result presents the main properties of quasi-Fej\'er convergent sequences; see  \cite[Theorem~$2.6$]{Bauschke:1996}.

\begin{lemma}[Quasi-F\'ejer Convergence Properties]\label{lem:quasi-Fejer} If
$(x_{k})_{k\in\mathbbm{N}}$ is quasi-Fej\'er convergent to $S$, then the following statements hold:
\item [ {\rm (a)}] The sequence $(x_{k})_{k\in\mathbbm{N}}$ is bounded;
\item [ {\rm (b)}] If an accumulation point $\bar x$ of $(x_{k})_{k\in\mathbbm{N}}$ belongs to $S$,
then $(x_{k})_{k\in\mathbbm{N}}$ is convergent to $\bar x$.
\end{lemma}
We conclude the section with a basic inequality that will be used in the next sections.
\begin{lemma}\label{ineq-norm} For any $x, y\in \mathbbm{E}$, we have
$\|x+y\|^2\le 2\|x\|^2+2\|y\|^2$.
\end{lemma}

\section{Inexact proximal solutions}\label{Inexact-Rules}

In this section, we introduce a concept of approximate solutions for the proximal gradient subproblem~\eqref{bar-xk-PGM}. We then describe how certain well-known algorithms compute these approximate solutions when the objective function exhibits a specific structure. Finally, we discuss the most closely related inexact proximal criteria proposed in the literature.

In the following, we introduce our concept of approximate proximal solution. First, recall that, given $x \in \mathbbm{E}$, the exact solution of subproblem  \eqref{bar-xk-PGM} with $\lambda_{k} = 1$ and $x_k=x$ consists of finding $\tilde{x}$ such that 
\begin{equation}\label{prox098}
 \tilde{x} = \argmin_{ y \in \mathbbm{E}} \left\{\left\langle \nabla f(x), y - x \right\rangle + g(y) + \frac{1}{2}\|y - x\|^2\right\}.  
\end{equation}
Equivalently, this corresponds to solving the following monotone inclusion problem:
\begin{equation}\label{incl:exactProxGrad-inclusion}
    0 \in \nabla f(x) + \partial g(\tilde{x}) + \tilde{x} - x.
\end{equation}

The concept of approximate solution given below consists of  relaxing the above  inclusion by introducing a residual pair \( (v, \varepsilon) \) that satisfies a specific mixed-relative error condition.

\begin{definition}\label{def:approx-prox-solution}
Let \( \tau \in (0,1] \), \( \gamma_1 > 1 \), \( \gamma_2 \geq 1 \), and \( \alpha \in [0, 1-\tau] \). Given a point \( x \in \mathbb{E} \), we say that \( \tilde{x} \in \mathbb{E} \) is an inexact proximal solution of  \eqref{incl:exactProxGrad-inclusion} if there exists a residual pair \( (v, \varepsilon) \in \mathbb{E} \times \RR_+ \) such that the triple \( (\tilde{x}, v, \varepsilon) \) satisfies the following conditions:
\begin{align}
&v \in \nabla f(x)+ \partial_{\ve}g( \tilde x) + \tilde x - x, \label{IR-inclusion}
\\
g( \tilde x-  v)-g( \tilde x)-\langle\nabla f(x),  v\rangle&+\frac{(1+\gamma_1)}{2}\|  v\|^2 +(1+\gamma_2)\ve\leq \frac{(1-\tau-\alpha)}{2}\|x-  \tilde x\|^2. \label{IR-inequality}
\end{align}
\end{definition}
\medskip

Let us now make  some remarks about the above definition. First note that the inclusion in \eqref{IR-inclusion} relaxes \eqref{incl:exactProxGrad-inclusion} by introducing an approximate solution $\tilde x$ together with a residues pair $(v,\varepsilon)$. The inequality in \eqref{IR-inequality} provides a mechanism for controlling the residual pair \((v, \varepsilon)\), the functional value of \(g\) at the approximate solution \(\tilde{x}\), and the angle between \(\nabla f(x)\) and the residual \(v\), all in terms of the distance between \(\tilde{x}\) and the point \(x\). Additionally,  \eqref{IR-inclusion}--\eqref{IR-inequality} ensure that both \(\tilde{x}\) and \(\tilde{x} - v\) lie within \(\dom g\). 

Now note that if $\tilde x$ is the exact solution of  the  monotone inclusion problem~\eqref{incl:exactProxGrad-inclusion}, we immediately have
\[
(\tilde x, v,\varepsilon):=\big(
\prox_{ g}
(x- \nabla f(x)),0, 0 \big)
\]
satisfies \eqref{IR-inclusion}-\eqref{IR-inequality}, for all $\tau \in  (0,1]$,  $\gamma_1>1$, $\gamma_2\ge 1$ and $\alpha\in [0,1-\tau]$.
Note also that if \( x \) is not an exact solution of \eqref{main-prob}, it cannot be the exact solution of \eqref{incl:exactProxGrad-inclusion} either.
Hence, by considering a point \( \tilde{x} \) that differs from $x$ but remains sufficiently close to the exact solution of \eqref{incl:exactProxGrad-inclusion}, it is direct that the right-hand side of \eqref{IR-inequality} is strictly positive because  \( 1 - \tau - \alpha > 0 \). Consequently, this inequality will eventually be satisfied by a residual pair \( (v, \varepsilon) \) that approaches zero, since its left-hand side also converges to zero in this case. 

It is worth pointing out that depending on the structure of the problem and/or the algorithm utilized for computing an inexact solution of the proximal subproblem,   the residual $v$ as in Definition~\ref{def:approx-prox-solution}  can be set as zero. In this case, \eqref{IR-inclusion}-\eqref{IR-inequality} are equivalent to 
\begin{gather}\label{eq:456}
0 \in \nabla f(x) + \partial_{\varepsilon} g(\tilde{x}) + \tilde{x} - x, \\
\varepsilon \leq \frac{(1 - \tau - \alpha)}{2(1 + \gamma_2)} \|x - \tilde{x}\|^2.\label{eq:4562}
\end{gather}

In the following two subsections, we focus on a specific structure of the nonsmooth component function  \( g \)  and provide a detailed discussion on how to compute a triple \( (\tilde{x}, v, \varepsilon) \) as described in Definition~\ref{def:approx-prox-solution}. Specifically, we address two scenarios: (i) \( g \) has a separable structure, which includes applications such as CUR-like factorization problems; and (ii) \( g \) is the indicator function of a special nonempty  convex and compact set \( C \).
We also note that inexact proximal solutions for the \( H \)-weighted nearest correlation matrix problem and the convex regularized problem can be obtained from \cite[Section~3]{Bello-Cruz:2023} through straightforward adaptations of the IR rules discussed therein.

\subsection{Inexact proximal solutions for separable functions}\label{sep.struct}

In this subsection, we discuss how inexact proximal solutions of \eqref{incl:exactProxGrad-inclusion} as in Definition~\ref{def:approx-prox-solution} can be computed if the nonsmooth component function $g$ has a separable structure of the form 
$
g = g_1 + g_2, 
$
where both $\prox_{g_1}$ and $\prox_{g_2}$ admit closed-form solutions.
 In this case, the proximal subproblem corresponds to  
\begin{equation}\label{prox09845}
\min_{y \in \mathbbm{E}} \left\{ g_1(y) + g_2(y) + \frac{1}{2}\|y - z\|^2 \right\},  
\end{equation}
where \( z = x - \nabla f(x) \).

We next discuss how an inexact proximal solution to \eqref{prox09845} can be obtained using either the Dykstra-like algorithm \cite[Theorem 3.3]{Bauschke:2008c} or the alternating direction method of multipliers (ADMM) \cite{Alves:2020, Adona:2019, Adona:2020}. { 
The approach related to the Dykstra algorithm is inspired by
\cite{Schmidt:2011} (see also \cite{Barre:2022, Zhou:2022}).
}

Let us first discuss the Dykstra-like algorithm. 
This algorithm, applied to problem \eqref{prox09845} with initial points \( z_0 = z \), \( p_0 = 0 \), and \( q_0 = 0 \), generates the sequences:  
\begin{equation}\label{eq:prox901}
\begin{aligned}
    y_\ell &= {\prox}_{g_1}(z_\ell + p_\ell), \quad &p_{\ell+1} &= z_\ell + p_\ell - y_\ell, \\ 
    z_{\ell+1} &= {\prox}_{g_2}(y_\ell + q_\ell), \quad &q_{\ell+1} &= y_\ell + q_\ell - z_{\ell+1}.
\end{aligned}
\end{equation}  
From the definition of \( y_\ell \) in \eqref{eq:prox901}, it follows that \( z_\ell + p_\ell - y_\ell \in \partial g_1(y_\ell) \).  
This observation, combined with Proposition~\ref{transp}, yields
\begin{equation}\label{eps}
z_\ell+p_\ell-y_\ell\in \partial_{\epsilon_\ell}  g_1(z_{\ell+1}), \quad \mbox{where} \quad 
\epsilon_\ell=g_1(z_{\ell+1})-g_1(y_{\ell})-\inner{
 z_\ell+p_\ell-y_\ell}{z_{\ell+1}-y_\ell}\geq 0.
\end{equation}
On the other hand, it follows from the definition of $z_{\ell+1}$ in \eqref{eq:prox901} that 
\(
y_\ell+q_\ell-z_{\ell+1}\in \partial g_2(z_{\ell+1}).
\)
Thus, by adding the inclusion in \eqref{eps} and the latter one, we have
\[
z_\ell+p_\ell+ q_\ell-  z_{\ell+1}\in \partial_{\epsilon_\ell} g(z_{\ell+1}).
\]
Rearranging the terms and using that $z=x-\nabla f(x)$, we get
\[
z_\ell+p_\ell+q_\ell-z\in \partial_{\epsilon_\ell} g(z_{\ell+1}) +z_{\ell+1}-x+\nabla f(x).
\]
It follows from the definitions of $p_{\ell+1}$ and $q_{\ell+1}$ in \eqref{eq:prox901} that $p_{\ell+1}+q_{\ell+1}+z_{\ell+1}=p_{\ell}+q_{\ell}+z_{\ell}$ for all $\ell$. Then, from the latter inclusion we have\[
z_0+p_0+q_0-z\in \partial_{\epsilon_\ell} g(z_{\ell+1}) +z_{\ell+1}-x+\nabla f(x).
\]
Since $z=p_{0}+q_{0}+z_{0}$, the inclusion becomes  
\[
0\in \partial_{\epsilon_\ell} g(z_{\ell+1}) +z_{\ell+1}-x+\nabla f(x).
\]
Finally, if the Dykstra-like algorithm is terminated when  
\begin{equation}\label{eq:o90}
   \varepsilon_\ell \leq \frac{(1-\tau-\alpha)}{2 (1+\gamma_2)} \|x - z_\ell\|^2, 
\end{equation}
then $(\tilde{x}, v, \varepsilon):=(z_\ell,0, \varepsilon_\ell) $   satisfies the conditions in \eqref{IR-inclusion}--\eqref{IR-inequality}. 

We also note that, as shown in \cite[Theorem 3.3]{Bauschke:2008c}, the sequences \( (y_{\ell})_{\ell \in\mathbbm{N}} \) and \( (z_{\ell})_{\ell \in\mathbbm{N}} \) converge to the exact solution \( \bar{x} \) of \eqref{prox09845}. Consequently, the sequence \( (\varepsilon_\ell)_{\ell\in\mathbbm{N}} \)  converges to zero provided that the subdifferential of \( g_1 \) is bounded over a compact set, which holds for instance if $\dom g_1$ is an open set. Furthermore, as long as \( \|x - \bar{x}\| > 0 \), it is guaranteed that the Dykstra-like algorithm will find a triple \( (\tilde{x}, 0, \varepsilon) \) that satisfies the conditions in \eqref{IR-inclusion}--\eqref{IR-inequality}, since the left-hand side of \eqref{eq:o90} will converge to zero
while the right-hand side remains strictly positive if $1-\tau-\alpha>0$.

Let us now discuss the ADMM. Note that \eqref{prox09845} can be rewritten as 
\begin{equation}\label{prox0984556}
\min \left\{ \hat g_1(x) + \hat g_2(y): x-y=0  \right\},  
\end{equation}
where \( \hat g_1(x) = g_1(x) \) and \( \hat g_2(y) = g_2(y) + \frac{1}{2}\|y - z\|^2 \). 
By applying the standard ADMM to the above problem, with initial point \( (y_0, \lambda_0) \in \mathbbm{E} \times \mathbbm{E} \) and penalty parameter \( c > 0 \), the following sequences are generated:
\begin{equation}\label{eq:prox90190}
    x_\ell = \prox_{g_1/c}\left(y_{\ell-1} + \frac{\lambda_{\ell-1}}{c}\right), \quad  
    y_{\ell} = \prox_{g_2/(c+1)}\left(\frac{z + c x_{\ell} - \lambda_{\ell-1}}{c+1}\right), \quad 
    \lambda_{\ell} = \lambda_{\ell-1} - c(x_\ell - y_\ell).
\end{equation}
From the definition of \( x_\ell \) in \eqref{eq:prox90190}, it follows that \( c y_{\ell-1} + \lambda_{\ell-1} - c x_\ell \in \partial g_1(x_\ell) \).  Combining the last inclusion with Proposition~\ref{transp}, we obtain  
\[
 c y_{\ell-1} + \lambda_{\ell-1} - c x_\ell \in \partial_{\varepsilon_\ell}  g_1(y_\ell), \quad \text{where} \quad 
\varepsilon_\ell = g_1(y_\ell) - g_1(x_{\ell}) - \langle  c y_{\ell-1} + \lambda_{\ell-1} - c x_\ell, y_\ell - x_{\ell} \rangle \geq 0.
\]
Similarly, from the definition of \( y_{\ell} \) in \eqref{eq:prox90190}, we have \( z + c x_{\ell} - \lambda_{\ell-1} - (c+1) y_{\ell} \in \partial g_2(y_{\ell}) \). 
Consequently, from the above inclusions, it follows that  
\[
c (y_{\ell-1}  -  y_{\ell})- y_\ell +z\in \partial_{\varepsilon_\ell} g(y_\ell),
\] 
which, combined with the definition of \( z \), yields 
\[
v_\ell := c (y_{\ell-1}  -  y_{\ell}) \in \partial_{\varepsilon_\ell} g(y_\ell) + y_\ell - x + \nabla f(x).
\] 
Finally, if the ADMM is terminated when  
\begin{equation}\label{eq:o901}
  g(y_\ell - v_\ell) - g(y_\ell) - \langle \nabla f(x), v_\ell \rangle + \frac{(1+\gamma_1)}{2} \|v_\ell\|^2 + (1+\gamma_2)\varepsilon_\ell \leq \frac{(1-\tau-\alpha)}{2} \|x - y_\ell\|^2,  
\end{equation} 
then $(\tilde{x}, v, \varepsilon):=(y_\ell,v_\ell, \varepsilon_\ell) $   satisfies the conditions in \eqref{IR-inclusion}--\eqref{IR-inequality}.

Note that, since the sequence \( ((x_{\ell}, y_{\ell}, \lambda_{\ell}))_{\ell \in\mathbbm{N}} \)  converges to  \( (\bar x, \bar y, \bar \lambda) \) which satisfies  \(0 \in \partial g_1(\bar x)-\bar \lambda \),  \(0 \in \partial g_2(\bar y)+\bar y-z-\bar \lambda  \) and $\bar x=\bar y$, it follows  that
 \( (v_\ell)_{\ell \in\mathbbm{N}} \) converges to zero, and \( (\varepsilon_\ell)_{\ell\in\mathbbm{N}} \) also converges to zero, as long as the subdifferential of \( g_1 \) is bounded over a compact set. Furthermore, as long as \( \|x - \bar{x}\| > 0 \), the ADMM is guaranteed to reach a triple \( (\tilde{x}, v, \varepsilon) \) satisfying the conditions in \eqref{IR-inclusion}--\eqref{IR-inequality}, since the left-hand side of \eqref{eq:o901} will converge to zero
while the right-hand side remains strictly positive provided that $1-\tau-\alpha>0$.

\subsection{Inexact projection  computed by the Frank-Wolfe method}
In this subsection, we assume that the nonsmooth component function $g$ is the indicator function of a nonempty,   convex and compact  set $C$. In this context, the proximal subproblem \eqref{prox098} (or equivalently the proximal monotone inclusion \eqref{incl:exactProxGrad-inclusion}) corresponds to project the point $x-\nabla f(x)$ onto the set $C$, and hence, the inexact proximal solution as in Definition~\ref{def:approx-prox-solution} can be seen as an inexact projection of  $x-\nabla f(x)$ onto the set $C$. By considering the case in which the residual $v$ can be set as $v=0$ (see \eqref{eq:456}-\eqref{eq:4562}), we  have that  \( \tilde{x} \in C \) is an inexact projection  with residual \( (0, \varepsilon) \in \mathbb{E} \times \RR_+ \) if and only if  
\begin{equation}\label{eq:0987}
\langle x - \nabla f(x) - \tilde{x}, y - \tilde{x} \rangle \leq \varepsilon\leq \frac{(1 - \tau - \alpha)}{2(1 + \gamma_2)} \|x - \tilde{x}\|^2, \quad \forall y \in C.
\end{equation}
For example,  one can use
the conditional gradient (CondG) method, a.k.a. Frank-Wolfe method \cite{Frank:1956,Jaggi:2013},
to compute $\tilde x \in C $ as in \eqref{eq:0987}. Indeed, 
given $z_{\ell}\in C$,  the $\ell$-th step of the CondG method,  applied to solve the projection problem
\begin{equation}\label{er67}
  \min_{y\in C}\frac{1}{2}\|y-x+  \nabla f(x)\|^2,  
\end{equation} first finds  $w_\ell$ as the minimizer of the linear function
$\langle{z_{\ell} - x+  \nabla f(x)},{\cdot - z_{\ell}} \rangle$ over  $ C$ and then sets $z_{\ell+1} =(1-\alpha_{\ell})z_{\ell}+\alpha_\ell w_\ell $, where $\alpha_\ell:=\min\{1,\langle{x-  \nabla f(x)-z_{\ell} },{w_{\ell} - z_{\ell}}\rangle/\|z_{\ell}-w_{\ell}\|^2\}$.
If, at iteration \( \ell \), the point \( z_\ell \) and  \( w_\ell \) satisfy the stopping criterion:  
\begin{equation}\label{eq:o9013}
\varepsilon_\ell :=\langle  x - \nabla f(x)-z_\ell, w_\ell - z_\ell \rangle \leq  \frac{(1 - \tau - \alpha)}{2(1 + \gamma_2)} \|x - z_\ell\|^2,
\end{equation} 
then \( \tilde{x} := z_\ell \) is an inexact  projection onto $C$ with residual  
$
(v, \varepsilon) := (0, \varepsilon_\ell).  
$ 
Furthermore, 
since the objective function in \eqref{er67} is strongly convex, the sequences \( (z_{\ell})_{\ell \in\mathbbm{N}} \) and \( (w_{\ell})_{\ell \in\mathbbm{N}} \) converge to \( P_C(x-\nabla f(x)) \). Hence, as long as \( \|x - P_C(x-\nabla f(x))\| > 0 \) and   \( 1-\tau-\alpha > 0 \), the conditional gradient method is guaranteed to reach a triple \( (\tilde{x}, 0, \varepsilon) \) satisfying the conditions in \eqref{eq:456}--\eqref{eq:4562}.  This is ensured because the left-hand side of \eqref{eq:o9013} converges to zero, while the right-hand side remains strictly positive.

\subsection{Related concepts of inexact proximal solutions}\label{HPE-inexact-subsection}
In this subsection, we discuss some works that analyze relative inexact proximal solution criteria, which are related to the concept introduced in Definition~\ref{def:approx-prox-solution}. Specifically, we review the criteria proposed in \cite{Solodov:1999a,Millan:2019,Bonettini:2016}.

We start by recalling the error criterion introduced in \cite{Solodov:1999a}, as applied to the proximal subproblem~\eqref{section1_1.1}. It is worth noting that this criterion was originally proposed in the broader context of monotone inclusions.  
 For a  given   point $x\in \mathbb{E}$, the proximal subproblem with a prox-stepsize $\lambda >0$ consists of computing $\tilde x$ such that
$$ 
0\in \lambda \partial g(\tilde x)+ \tilde x-x.
$$
The criterion  in \cite{Solodov:1999a} relaxes this inclusion by finding a triple \((\tilde{x}, w, \varepsilon) \in \mathbb{E} \times \mathbb{E} \times \RR_+\) that satisfies  
\begin{equation}\label{hpe-conditions}
w \in \partial_\varepsilon g(\tilde{x}), \qquad \|\lambda w + \tilde{x} - x\|^2 + 2\lambda\varepsilon \leq \sigma\|\tilde{x} - x\|^2,  
\end{equation}  
for some scalar \(\sigma \in [0,1)\).

The hybrid proximal extragradient (HPE) method proposed in \cite{Solodov:1999a} consists mainly of two basic steps: first, it computes an approximate proximal solution based on  \eqref{hpe-conditions}, and then an extragradient step is performed to determine the next iterate, given by \(x^+ = x - \lambda w\). The convergence properties of the HPE method were analyzed in \cite{Solodov:1999a}, while its iteration complexity was established in \cite{Monteiro:2010}.
As demonstrated in \cite{Solodov:1999a,Monteiro:2010}, the HPE method can be seen as a general framework that encompasses various well-known algorithms, including the Korpelevich method \cite{Korpelevich:1976} and Tseng's modified forward-backward splitting (Tseng-MFBS) method \cite{Tseng:2000}, among others.

Tseng-MFBS method  can be, in particular, applied for solving  \eqref{main-prob}. In this case, it  generates a sequence \((x_k)_{k\in\mathbbm{N}}\) as follows:  
\begin{equation}\label{Tseng-method}
x_{k+1} = \bar{x}_k - \lambda (\nabla f(\bar{x}_k) - \nabla f(x_k)), \qquad \bar{x}_k := \prox_{\lambda g}(x_k - \lambda \nabla f(x_k)), 
\end{equation}  
where $x_0$ is a given initial point, $\lambda > 0$ is the stepsize, and $\prox_{\lambda g}$ is defined in \eqref{section1_1.1}.  

As shown in \cite{Monteiro:2010,Solodov:1999a}, this scheme is an instance of the HPE method when $\nabla f$ is $L$-Lipschitz continuous. Indeed, by defining  
\[
w_k := \frac{1}{\lambda}(x_k - \bar{x}_k) + \nabla f(\bar{x}_k) - \nabla f(x_k),
\]  
it can be verified that $(\tilde{x}, w, \varepsilon) := (\bar{x}_k, w_k, 0)$ satisfies the conditions in \eqref{hpe-conditions} with $\sigma = \lambda / L$.

Most recently, \cite{Millan:2019} proposed   inexact  PGMs for solving  \eqref{main-prob},  where both component functions \(f\) and \(g\) may be nonsmooth.  
Algorithm~2 in this reference (with a fixed stepsize \(\alpha_k = \lambda\) and  \(C = \mathbb{E}\)) applied to \eqref{main-prob}, with the additional restriction that $\dom g=\mathbb{E}$, generates a sequence \((x_k)_{k\in\mathbbm{N}}\) according to the following scheme: given \(x_k \in \mathbb{E}\) and \(\sigma \in [0,1)\), first let \(z_k := x_k - \lambda \nabla f(x_k)\), and then compute a triple \((\tilde{x}_k, w_k, \varepsilon_k) \in \mathbb{E} \times \mathbb{E} \times \RR_+\) such that
$$
w_k \in \partial_{\varepsilon_k} g(\tilde{x}_k), \qquad \|\lambda w_k + \tilde{x}_k - z_k\|^2 + 2\lambda \varepsilon_k \leq \sigma \|\tilde{x}_k - z_k\|^2,
$$
and update \(x_k\) as \(x_{k+1} = z_k - \lambda w_k\).  
Note that, by defining $v_k = \lambda w_k + \tilde{x}_k - z_k$, the above  conditions, combined with the definition of $z_k$, are equivalent to  
$$
\frac{1}{\lambda} v_k \in \nabla f(x_k) + \partial_{\varepsilon_k} g(\tilde{x}_k) + \frac{1}{\lambda} (\tilde{x}_k - x_k), \qquad \|v_k\|^2 + 2\lambda \varepsilon_k \leq \sigma^2 \|\tilde{x}_k - z_k\|^2.
$$  
Note that, the above inclusion corresponds to the one in \eqref{IR-inclusion} if $\lambda = 1$. On the other hand,  while the latter inequality represents a relative inexact proximal criterion, it differs from the one in \eqref{IR-inequality}.
It is worth noting that \cite{Millan:2019} did not incorporate any linesearch procedure.
Instead, they considered stepsizes such as fixed, exogenous, and Polyak-type.

We end this subsection by discussing \cite{Bonettini:2016}, which also studies a PGM with inexact proximal solutions and an explicit linesearch. To simplify the discussion, we assume that the prox-parameter in \cite{Bonettini:2016} is constant and equal to one, and that the variable metric, used to define their proximal subproblem, is fixed and coincides with the Euclidean one.
In this case, at the $k$-th iteration, they compute a pair \((\tilde{x}_k,  \varepsilon_k) \in \mathbb{E} \times  \RR_+\) such that 
\begin{equation}\label{inclusion-with-0}
0 \in \nabla f(x_k) + \partial_{\varepsilon_k} g(\tilde{x}_k) + \tilde{x}_k - x_k,
\end{equation}  
i.e., the triple $(\tilde{x}_k, 0, \varepsilon_k)$ is a solution of \eqref{IR-inclusion} (with a null residual $v_k$). 
Two types of conditions are assumed to control the residual sequence $(\varepsilon_k)_{k\in\mathbbm{N}}$: (i) the summable absolute error condition; and (ii) a relative error condition. We focus our attention on the second one, which is more closely related to our inexact criterion. The authors assume that the pair $(\tilde{x}_k, \varepsilon_k)$ satisfies \eqref{inclusion-with-0} and the following conditions:  
\begin{equation}\label{prox-condition-VarMetrPaper}
\varepsilon_k \leq -t \tilde{h}_s(\tilde{x}_k), \qquad \tilde{h}_s(\tilde{x}_k) < 0,
\end{equation}  
for some $s \in (0, 1]$ and $t > 0$, where  
\begin{equation}\label{h_s}
\tilde{h}_s(u) = \langle \nabla f(x_k), u - x_k \rangle + g(u) - g(x_k) + \frac{s}{2} \|u - x_k\|^2.
\end{equation}
Note that $\tilde{h}_s$ is closely related to the objective function of the subproblem \eqref{prox098} and becomes negative at its unique solution. In particular, the second condition in \eqref{prox-condition-VarMetrPaper} holds for some $\tilde{x}_k$ sufficiently close to the exact solution of \eqref{prox098}. 

We claim that if $(\tilde{x}, v, \varepsilon) := (\tilde{x}_k, 0, \varepsilon_k)$ is an approximate proximal solution in the sense of Definition~\ref{def:approx-prox-solution}, then the pair $(\tilde{x}_k, \varepsilon_k)$ satisfies \eqref{inclusion-with-0}-\eqref{prox-condition-VarMetrPaper} for some $t>0$. Indeed, in this case,  \eqref{IR-inclusion} becomes 
$$
0\in \nabla f(x)+\partial_{\varepsilon} g(\tilde x)+\tilde x-x
$$
and \eqref{IR-inequality}  is equivalent to 
\begin{equation}\label{eps-ineq-nullv}
\ve\leq \frac{1-\tau-\alpha}{2(1+\gamma_2)}\|x- \tilde x\|^2.
\end{equation}
The above inclusion, together with the definitions of $\tilde h_s$ in \eqref{h_s} and the   $\varepsilon$-subdifferential in \eqref{def-epsSubdiff},  imply that
$$
\tilde h_s(\tilde x)+\left(\frac{2-s}{2}\right)\|x-\tilde x\|^2\leq \ve.
$$
Combining the latter two inequalities,  we obtain
\begin{equation*}
\tilde h_s(\tilde x)\leq \left[\frac{1-\tau-\alpha}{2(1+\gamma_2)}-\frac{2-s}{2}\right]\|x-\tilde x\|^2.
\end{equation*}
Hence, if $t$ is a scalar such that
\begin{equation}\label{ineq-t}
t\geq \frac{1-\tau-\alpha}{(2-s)(1+\gamma_2)+\tau+\alpha-1}.
\end{equation}
Thus, in view of  the latter inequality and \eqref{eps-ineq-nullv} that 
$$
  \varepsilon\leq \frac{1-\tau-\alpha}{2(1+\gamma_2)}\|x- \tilde x\|^2\leq -t\tilde h_s(\tilde x),
$$
which implies that the pair $(\tilde x,\varepsilon)$ satisfies the conditions in \eqref{prox-condition-VarMetrPaper} with $t$ as in \eqref{ineq-t}, proving the claim. 

We conclude by observing that even if the residue $v$ in Definition \ref{def:approx-prox-solution} is chosen as null, then the IPG-ELS algorithm differs from the scheme proposed in\cite{Bonettini:2016}. This is because the linesearch we propose is applied only to the smooth part $f$ of the objective function $F$, whereas the Armijo-type linesearch procedure considered in \cite{Bonettini:2016} is applied to the whole objective function $F$. Moreover, if $v \neq 0$, then our inexact criterion differs significantly from the one in \cite{Bonettini:2016}. In addition, our linesearch procedure accounts for not only the proximal inexact solution $\tilde x$ but also the residue pair $(v, \varepsilon)$. It is also worth mentioning that the latter reference did not analyze the iteration complexity to achieve an approximate stationary solution in the sense of Definition~\ref{approxStationarySolution}. Finally, \cite{Bonettini:2016} also proposes a second criterion, where an approximate proximal solution $\tilde x$ of the subproblem \eqref{incl:exactProxGrad-inclusion} is accepted if
$$
\tilde h_1(\tilde x) \leq \eta \tilde h_1(\bar x),
$$ where $\eta\in(0,1]$ 
and $\bar x$ being is exact solution of \eqref{incl:exactProxGrad-inclusion}.

\section{Inexact Proximal Gradient Method}\label{search}
In this section, we introduce the inexact proximal gradient method with an explicit linesearch and establish some basic properties.

In the following, we formally describe the proposed method.
\begin{center}\fbox{\begin{minipage}[b]{\textwidth}
{\bf Inexact Proximal Gradient Method with an Explicit  Linesearch (IPG-ELS) }
\item [\;1.] {\bf Initialization Step.} Let $x_0\in \dom g$, $\tau \in (0,1]$, $\theta \in (0,1)$,   $\gamma_1>1$, $\gamma_2\ge 1$, and  $\alpha\in [0, 1-\tau]$. Set $k=0$.
\item [\;2.] {\bf Inexact Proximal Subproblem.}  {Compute an inexact proximal solution \( \tilde{x}_k \) along with a residual pair \( (v_k, \varepsilon_k) \) for \eqref{incl:exactProxGrad-inclusion} with \( x := x_k \). Specifically, find a triple \( (\tilde{x}_k, v_k, \varepsilon_k) \) such that}
\begin{gather}\label{inexProx-inclusion}
v_{k}\in \nabla f(x_{k})+\partial_{\varepsilon_k}g( \tilde x_{k})+\tilde x_{k}- x_{k}, \\[2mm]\label{inex-relative-condition}
g(\tilde x_{k}-  v_{k})-g(\tilde x_{k})-\langle\nabla f(x_{k}),  v_{k}\rangle+\frac{(1+\gamma_1)}{2}\|  v_{k}\|^2 +(1+\gamma_2)\ve_{k}\leq \frac{(1-\tau-\alpha)}{2}\|x_{k}- \tilde x_{k}\|^2;
\end{gather}
\item [\;3.] \noindent  {\bf Stopping Criterion.}  If $\tilde x_{k}=x_{k}$, then stop.
\item [\;4.] \noindent {\bf Linesearch Procedure.}  Set $\beta=1$ and $y_{k}:=\tilde x_{k}-  v_{k}$.
  If 
\begin{equation}\label{ineq:linesearch}\dsty f\left(x_{k} +\beta (y_{k}-x_{k})\right) \leq f(x_{k})+\beta
\la\nabla f(x_{k}), y_{k}-x_{k}\ra+\frac{\beta\tau}{2}\|x_{k}-\tilde x_{k}\|^2+\frac{\gamma_1}{2}\beta\|v_{k}\|^2+\beta\gamma_2 \varepsilon_k,
\end{equation}
then set $\beta_{k}=\beta$, $x_{k+1}=x_{k} +\beta_{k} (y_{k}-x_{k})$ and $k:=k+1$, and go to Step~2. Otherwise, set $\beta=\theta \beta$, and verify~\eqref{ineq:linesearch}.
\end{minipage}}\end{center}

\begin{remark}[The Explicit Lineasearch]
Note that the novel explicit linesearch in Step 4 of the IPG-ELS is related to the ones proposed in \cite{Salzo:2017a,BelloCruz:2016a} when \(v_{k} = 0\) and \(\varepsilon_k = 0\). Moreover, note that the linesearch does not evaluate the inexact proximal subproblem inside  the inner loop of Step 4. Hence, only one inexact proximal solution is computed per iteration.  
\end{remark}
 
We introduce a useful lemma that plays a crucial role in analyzing the stopping criterion of the IPG-ELS method.

\begin{lemma}[Iteration Inequality Condition] The following inequality holds for every iteration $k$.
\begin{equation}\label{ineq-useful-lem-stop}
\left\langle \tilde x_{k}- x_{k},v_{k}\right\rangle+\frac{(\gamma_1-1)}{2}\| v_{k}\|^2+\gamma_2 \varepsilon_k\leq \frac{(1-\tau-\alpha)}{2}\|x_{k}- \tilde x_{k}\|^2.
\end{equation}
  \end{lemma}  
{\it Proof} It follows from \eqref{inexProx-inclusion} that 
\[
v_{k}+ {x}_{k}-\tilde x_{k}-\nabla f({x}_{k})\in \partial_{\varepsilon_k}g( \tilde x_{k}).\]
Now using
\eqref{def-epsSubdiff} together with the fact that $y_{k}=\tilde x_{k}-v_{k}$, we have  
\begin{align*}
g(y_{k})-g( \tilde x_{k})&\ge \left \langle v_{k}+ {x}_{k}-\tilde x_{k}-\nabla f({x}_{k}),y_{k}- \tilde x_{k}\right \rangle-\ve_{k}\\&=\left \langle v_{k}+ {x}_{k}-\tilde x_{k}-\nabla f({x}_{k}),-v_{k}\right \rangle-\ve_{k},
\end{align*}
which is equivalent to 
\begin{align}\label{lkge}
\Gamma_{k}:=g(y_{k})-g( \tilde x_{k})-\langle \nabla f({x}_{k}),v_{k}\rangle
&\ge 
-\|v_{k}\|^2+\langle \tilde{x}_{k}-  x_{k},v_{k}\rangle-\ve_{k}.
\end{align}
On the other hand, considering the definitions of $\Gamma_{k}$ and $y_{k}$, we observe that \eqref{inex-relative-condition} is equivalent to 
\[
\Gamma_{k} +\frac{(1+\gamma_1)}{2}\| v_{k}\|^2+(1+\gamma_2) \varepsilon_k\leq \frac{(1-\tau-\alpha)}{2}\|x_{k}- \tilde x_{k}\|^2,
\]
which,  combined with \eqref{lkge}, yields the desired inequality \eqref{ineq-useful-lem-stop}.
\qed

The following result demonstrates that the termination criterion of the IPG-ELS method, as specified in Step 3, is satisfied only when a solution to problem~\eqref{main-prob} is identified.

\begin{lemma}[Stopping at a Solution]\label{lem:stopping at solution}
The IPG-ELS method terminates at the $k$-th iteration if and only if $x_{k}$ is a solution to problem~\eqref{main-prob}.
\end{lemma}
{\it Proof} 
Assume that the IPG-ELS method stops at the $k$-th iteration. In view of Step~3, we have $\tilde{x}_{k} = x_{k}$. Hence, it follows from \eqref{ineq-useful-lem-stop} that
\[
\frac{(\gamma_1 - 1)}{2}\|v_{k}\|^2 + \gamma_2 \varepsilon_k \leq 0.
\]
Since $\gamma_1 > 1$ and $\gamma_2 \geq 1$, we obtain $v_{k} = 0$ and  $\varepsilon_k = 0$. Hence, in view of \eqref{inexProx-inclusion}, we get that $0 \in \nabla f(x_{k}) + \partial g(x_{k}),$ concluding that $x_{k}$ is a solution of problem~\eqref{main-prob}.

Assume now that $x_{k}$ is a solution of problem~\eqref{main-prob}. Thus,
\(
-\nabla f(x_{k}) \in \partial g(x_{k}).
\)
It follows from \eqref{inexProx-inclusion} that $v_k+x_{k} - \tilde{x}_{k} - \nabla f(x_{k}) \in \partial_{\varepsilon_k} g(\tilde{x}_{k})$. So,  the $\varepsilon$-monotonicity of $\partial_\varepsilon g$ in \eqref{eps-monotonicity of subdiff g} implies
\[
\langle v_{k} + x_{k} - \tilde{x}_{k}, \tilde{x}_{k} - x_{k} \rangle \geq -\varepsilon_k,
\]
which is equivalent to
\begin{equation}
\langle v_{k}, \tilde{x}_{k} - x_{k} \rangle \geq \| x_{k}-\tilde{x}_{k} \|^2 - \varepsilon_k.
\end{equation}
Hence, it follows from \eqref{ineq-useful-lem-stop} that
\[
\frac{(1 + \tau + \alpha)}{2}\|x_{k} - \tilde{x}_{k}\|^2 + \frac{(\gamma_1 - 1)}{2}\|v_{k}\|^2 + (\gamma_2 - 1) \varepsilon_k \leq 0.
\]
Since $\gamma_1> 1$, $\gamma_2\ge 1$, $\alpha \geq 0$, and $\tau > 0$, we conclude that $x_{k} = \tilde{x}_{k}$ and $v_{k} = 0$. Hence, from \eqref{inex-relative-condition}, we also have $\varepsilon_k = 0$.  Therefore, the IPG-ELS method stops at the $k$-th iteration.
\qed

In the following, we establish the well-definedness of the linesearch procedure in Step~4 of the IPG-ELS method. 

\begin{lemma}[Finite Linesearch Termination]\label{boundary-well2} The linesearch procedure in Step~4 stops after a finite number of steps.
\end{lemma}
{\it Proof} 
In view of Step~3, we have $x_{k} \neq \tilde{x}_{k}$. Assume, for the sake of contradiction, that the linesearch procedure does not terminate after a finite number of steps. Thus, for all $\beta \in  \{1, \theta, \theta^2, \ldots\}$, we have
\begin{equation*}
f\left(x_{k} + \beta (y_{k} - x_{k})\right) > f(x_{k}) + \beta \langle \nabla f(x_{k}), y_{k} - x_{k} \rangle + \frac{\beta \tau }{2}\|x_{k} - \tilde{x}_{k}\|^2 + \frac{\gamma_1}{2}\beta\|v_{k}\|^2 + \beta \gamma_2 \varepsilon_k,
\end{equation*}
or, equivalently,
\[
\frac{f\left(x_{k} + \beta (y_{k} - x_{k})\right) - f(x_{k})}{\beta} - \langle \nabla f(x_{k}), y_{k} - x_{k} \rangle > \frac{\tau}{2}\|x_{k} - \tilde{x}_{k}\|^2 + \frac{\gamma_1}{2} \|v_{k}\|^2 + \gamma_2 \varepsilon_k.
\]
Given that $f$ is differentiable, the left-hand side of the above inequality approaches zero as $\beta \downarrow 0$, leading to the conclusion that
\begin{align*}
0 \geq \frac{\tau}{2}\|x_{k} - \tilde{x}_{k}\|^2 + \frac{\gamma_1}{2} \|v_{k}\|^2 + \gamma_2 \varepsilon_k.
\end{align*}
This implies $x_{k} = \tilde{x}_{k}$ contradicting the assumption that $x_{k} \neq \tilde{x}_{k}$.
\qed 

The above lemma ensures that the linesearch procedure in Step~4 of the IPG-ELS method terminates after a finite number of iterations. This result is fundamental to the convergence analysis of the method, as it guarantees that the linesearch procedure is well-defined.

Next, we provide a remark on the possibility of initiating the linesearch with a larger stepsize.

\begin{remark}[Over-relaxation Strategy] In the current implementation, the linesearch algorithm begins with an initial stepsize $\beta = 1$. Considering larger initial values for $\beta$, thereby adopting an over-relaxation strategy, may potentially accelerate convergence by extrapolating beyond standard update steps. This approach is related to the concept of conic averagedness \cite{BartzDaoPhan}, which offers a framework for analyzing the convergence properties of fixed-point algorithms under relaxed conditions. However, initiating the linesearch with $\beta > 1$ presents significant challenges to our current approach. While the linesearch remains well-defined for any $\beta > 1$, and the proof of finite termination would be directly extended, the convergence and complexity analysis of our scheme depend on the convexity of the function $f$ and the appropriateness of the direction $y_k - x_k$ to ensure a sufficient decrease in the objective function. Therefore, our convergence analysis does not cover such cases, and the practical implications for this approach are unclear. 
\end{remark}

The subsequent analysis investigates the complexity of the linesearch procedure introduced in Step~4 of the IPG-ELS method. Here we assume that the gradient of the function $f$, $\nabla f$, is Lipschitz continuous only to establish an upper bound for the number of iterations required by the linesearch procedure in Step~4.

\begin{lemma}[Lipschitz Condition and Linesearch Complexity]\label{lem:Linesearch-Lipschitz}
Assume that \(f\) has an \(L\)-Lipschitz continuous gradient and that \(x_{k}\) is not a solution to problem~\eqref{main-prob}. Then, any \(\beta \leq \tau/(2L)\) satisfies \eqref{ineq:linesearch}. As a consequence, the linesearch procedure in Step~4 of the IPG-ELS method stops in at most 
\begin{equation}\label{numb-iter-procedure}
\ell := \left\lceil \frac{\ln (\min\{\tau/(2L), 1\})}{\ln(\theta)} \right\rceil
\end{equation}
iterations.
\end{lemma}

{\it Proof} 
Since \(\nabla f\) is \(L\)-Lipschitz continuous, for any \(\beta > 0\), we have
\begin{align*}
f(x_{k} + \beta (y_{k} - x_{k})) - f(x_{k}) - \beta \langle \nabla f(x_{k}), y_{k} - x_{k} \rangle &\leq \frac{L \beta^2}{2}\|y_{k} - x_{k}\|^2.
\end{align*}
Hence, if \(\beta \leq \tau/(2L)\), we conclude that
\begin{align*}
f(x_{k} + \beta (y_{k} - x_{k})) - f(x_{k}) - \beta \langle \nabla f(x_{k}), y_{k} - x_{k} \rangle &\leq \frac{\tau \beta}{4}\|y_{k} - x_{k}\|^2= \frac{\tau \beta}{4}\|x_{k}-\tilde x_{k}+v_{k}\|^2\\&\le   \frac{\beta\tau}{2} \|x_{k}-\tilde x_{k}\|^2+\frac{\tau}{2} \beta \|v_{k}\|^2,
\end{align*} using Lemma \ref{ineq-norm} in the last inequality.
Since \(\gamma_1>1>\tau\) and  \(\gamma_2 \ge 0\), we have that  \eqref{ineq:linesearch} holds, thereby proving the first statement of the lemma. The last statement follows from the first one, given that the natural number \(\ell\), defined in \eqref{numb-iter-procedure}, satisfies \(\beta_{\ell} := \theta^\ell \leq \min\{\tau/(2L), 1\}\).
\qed 

This lemma provides a sufficient condition to ensure that the lower bound of the sequence generated by our linesearch is strictly greater than $0$. Specifically, if $\nabla f$ is $L$-Lipschitz continuous, then the stepsizes $\beta_k$, produced through the linesearch \eqref{ineq:linesearch}, are guaranteed to be bounded below by a positive constant $\beta > 0$, i.e., $\beta_k \ge \beta$ for all $k \in \mathbb{N}$. Moreover, it is possible to relax the global Lipschitz condition to something local, such as $\nabla f$ being locally Lipschitz continuous around any solution of problem~\eqref{main-prob}, as was done in Proposition 5.4(ii) of \cite{BelloCruz:2016a}. In fact, the proof of Lemma \ref{lem:Linesearch-Lipschitz} may be readily adapted to establish the same complexity now with respect to the locally Lipschitz constant $\mathcal{L}$.
\begin{lemma}[Locally Lipschitz Condition and Linesearch Complexity]\label{lem:Local-Linesearch-Lipschitz}
If \(\nabla f\) is \(\mathcal{L}\)-locally Lipschitz continuous at any solution of problem~\eqref{main-prob}, then any \(\beta \leq \tau/(2\mathcal{L})\) satisfies \eqref{ineq:linesearch}. Consequently, the linesearch procedure in Step~4 of the IPG-ELS method terminates after at most 
\begin{equation}\label{numb-iter-procedure}
\ell := \left\lceil \frac{\ln (\min\{\tau/(2\mathcal{L}), 1\})}{\ln(\theta)} \right\rceil
\end{equation}
iterations.
\end{lemma}
Note that the assumption of the gradient of $f$ in problem \eqref{main-prob} being locally, rather than globally, Lipschitz continuous is commonly encountered in practice. For instance, this condition arises in the Poisson linear inverse regularization problem  with Kullback–Leibler divergence \cite{Csiszar:1991,Vardi:1985}. We emphasize that this weaker assumption is required exclusively for analyzing the complexity of the linesearch procedure, as it ensures that the stepsizes generated by the proposed linesearch are uniformly bounded away from zero. Additionally, it may be used as sufficient condition for establishing the convergence rate of the functional value sequence. It is worth noting that the finite termination of the linesearch and convergence of the IPG-ELS method do not rely on such assumption.

\section{Convergence and Complexity Analysis of the IPG-ELS Method}\label{convSec}
In this section, we focus on analyzing the convergence properties of the IPG-ELS method. We establish the convergence and its iteration complexity in terms of functional values and a residual associated with an approximate solution, as defined in Definition~\ref{approxStationarySolution}.

We begin this section by presenting a result that is fundamental for establishing the convergence and the iteration complexity of the IPG-ELS method.

  \begin{proposition}[Key Inequality for the IPG-ELS Method]\label{Prop:Main-proposition}
For every $x\in \dom \, g$ and $ k\in \mathbbm{N}$, we have
\begin{equation}\label{ineq:MainIneq}
2\beta_{k}[F( x_{k})-F(x)]+  2({F(x_{k+1})-F(x_{k})})\leq   
 \|x_{k}-x\|^ 2-\|x_{k+1}-x\|^2- {\alpha\beta_{k}}\|x_{k}- \tilde x_{k}\|^2.
\end{equation}
Additionally, the sequence $(F( x_{k}))_{k\in\mathbbm{N}}$ is decreasing and convergent, and $\sum_{k=0}^{+\infty} \|x_k-x_{k+1}\|^2<+\infty$. 
\end{proposition} 
{\it Proof} 
Let $x\in \dom \, g$ and $k \in \mathbb{N}$. In view of the inexact relative inclusion \eqref{inexProx-inclusion}, we have
 $v_{k}+{x}_{k}- \tilde x_{k}-\nabla f({x}_{k})\in\partial_{\ve_{k}} g( \tilde x_{k})$, and hence the definition of $\partial_\ve g$ in \eqref{def-epsSubdiff} implies that
\begin{equation}\label{eq2-2}
g(x)\geq g( \tilde x_{k}) + \langle v_{k}+{x}_{k}- \tilde x_{k}-\nabla f({x}_{k}),x- \tilde x_{k}\rangle-\ve_{k}.
\end{equation}
Since $f$ is convex, we have 
$f(x)-f(x_{k})\ge \langle\nabla f(x_{k}), x-x_{k}\rangle$.
Adding the above two inequalities, using $f+g=F$, and simplifying the resulting expression, we obtain
$$
F(x)-F(x_{k}) \geq g(\tilde x_{k})- g(x_{k})+ \langle\nabla f(x_{k}), \tilde x_{k}-x_{k}\rangle - \inner{{x}_{k}- \tilde x_{k}}{ \tilde x_{k}-x}+\inner{v_{k}}{x-\tilde x_{k}}-\ve_{k}.
$$
Combining the above inequality with the identity
\[-\inner{{x}_{k}- \tilde x_{k}}{ \tilde x_{k}-x} =\frac{1}{2}\left[\|x_{k}- \tilde x_{k}\|^2 +\|\tilde x_{k}-x\|^2-\|x_{k}-x\|^2\right],\]
we have
\begin{align*}
F(x)-F( x_{k})  \geq &  \; g(\tilde x_{k})- g( x_{k})+
\langle\nabla f(x_{k}), \tilde x_{k}-x_{k}\rangle
+ \frac{1}{2}\|x_{k}-\tilde x_{k}\|^2\\
& + \frac{1}{2}\left[\|\tilde x_{k}-x\|^2-\|x_{k}-x\|^ 2\right]+\inner{v_{k}}{x- \tilde x_{k}}-\ve_{k},
\end{align*}
or, equivalently, 
\begin{align*}
F(x)-F( x_{k})+\frac{F(x_{k})-F(x_{k+1})}{\beta_{k}}\geq &\frac{g(x_{k})-g(x_{k+1})}{\beta_{k}}+g(\tilde x_{k})- g( x_{k})+ \langle\nabla f(x_{k}), \tilde x_{k}-y_{k}\rangle\\
&+\frac{f(x_{k})-f(x_{k+1})}{\beta_{k}}+ \langle\nabla f(x_{k}), y_{k}-x_{k}\rangle+ \frac{1}{2}\|x_{k}- \tilde x_{k}\|^2\\
&+ \frac{1}{2}\left[\|x  -(\tilde x_{k}- v_{k})\|^2-\|x_{k}-x\|^ 2\right] -\frac{1}{2}\|  v_{k}\|^2 -\ve_{k}.
\end{align*}
Since  $x_{k+1}=x_{k} +\beta_{k} (y_{k}-x_{k})$ and $g$ is convex, we have $g(x_{k+1})-g(x_{k})\leq\beta_{k}(g(y_{k})-g(x_{k}))$. Hence, combining the last two inequalities and the fact that $y_{k}=\tilde x_{k}-v_{k}$, we get
\begin{align*}
F(x)-F( x_{k})+\frac{F(x_{k})-F(x_{k+1})}{\beta_{k}}\geq & \frac{1}{2}\left[\|x-y_{k}\|^2-\|x_{k}-x\|^ 2\right]\\
&+\frac{f(x_{k})-f(x_{k+1})}{\beta_{k}}+ \langle\nabla f(x_{k}), y_{k}-x_{k}\rangle+ \frac{1}{2}\|x_{k}- \tilde x_{k}\|^2\\
&+g(\tilde x_{k})- g( y_{k})+ \langle\nabla f(x_{k}), \tilde x_{k}-y_{k}\rangle -\frac{1}{2}\|  v_{k}\|^2 -\ve_{k}.
\end{align*}
or, equivalently,

\begin{align*}
F(x)-F( x_{k})&+\frac{F(x_{k})-F(x_{k+1})}{\beta_{k}}\geq  \frac{1}{2}\left[\|x  -y_{k}\|^2-\|x_{k}-x\|^ 2\right]+ \frac{\alpha}{2}\|x_{k}- \tilde x_{k}\|^2\\
&+\frac{f(x_{k})-f(x_{k+1})}{\beta_{k}}+ \langle\nabla f(x_{k}), y_{k}-x_{k}\rangle+ \frac{\tau}{2}\|x_{k}- \tilde x_{k}\|^2+\frac{\gamma_1}{2}\|  v_{k}\|^2 +\gamma_2\ve_{k}\\
&+g(\tilde x_{k})- g( y_{k})+ \langle\nabla f(x_{k}), \tilde x_{k}-y_{k}\rangle -\frac{(1+\gamma_1)}{2}\|  v_{k}\|^2 -(1+\gamma_2)\ve_{k}\\&+ \frac{(1-\tau-\alpha)}{2}\|x_{k}- \tilde x_{k}\|^2.
\end{align*}
Now, using the linesearch procedure of Step~4, we obtain

\begin{align*}
F(x)-F( x_{k})+\frac{F(x_{k})-F(x_{k+1})}{\beta_{k}}\geq &  \frac{1}{2}\left[\|x  -y_{k}\|^2-\|x_{k}-x\|^ 2\right]+ \frac{\alpha}{2}\|x_{k}- \tilde x_{k}\|^2\\
&+g(\tilde x_{k})- g( y_{k})+ \langle\nabla f(x_{k}), \tilde x_{k}-y_{k}\rangle \\&-\frac{(1+\gamma_1)}{2}\|  v_{k}\|^2 -(1+\gamma_2)\ve_{k}+ \frac{(1-\tau-\alpha)}{2}\|x_{k}- \tilde x_{k}\|^2.
\end{align*}
It follows from the above inequality and \eqref{inex-relative-condition} that
\begin{equation}\label{ineq:aux-end of lem}
\beta_{k}[F(x)-F( x_{k})]+  {F(x_{k})-F(x_{k+1})}\geq   
 \frac{\beta_{k}}{2}\left[\|x  -y_{k}\|^2-\|x_{k}-x\|^ 2\right]+ \frac{\alpha\beta_{k}}{2}\|x_{k}- \tilde x_{k}\|^2.
\end{equation}
On the other hand, using the identity $x_{k+1}-x=(1-\beta_{k})(x_{k}-x) +\beta_{k} (y_{k}-x)$ and the strong convexity of $\|\cdot\|^2$, we have 
\[
\|x_{k+1}-x\|^2 \leq (1-\beta_{k})\|x_{k}-x\|^2 +\beta_{k} \|y_{k}-x\|^2-(1-\beta_{k})\beta_{k}\|x_{k}-y_{k}\|^2,
\]
which implies 
\[
\beta_{k}\left(\|y_{k}-x\|^2-\|x_{k}-x\|^2\right)\geq \|x_{k+1}-x\|^2-\|x_{k}-x\|^2.
\]
Therefore, the proof of \eqref{ineq:MainIneq}  follows by combining the latter inequality with  \eqref{ineq:aux-end of lem}. The last statement of the proposition follows immediately from \eqref{ineq:MainIneq}  with $x=x_{k}$ that \begin{equation}\label{ineq:MainIneq*}
2({F(x_{k+1})-F(x_{k})})\leq   
 -\|x_{k+1}-x_k\|^2<0.
\end{equation}
So, the sequence $(F(x_{k}))_{k\in\mathbbm{N}}$ is decreasing and convergent because it is bounded from below by $F(x_*)$. Moreover, the last inequality implies that $\sum_{k=0}^{+\infty} \|x_k-x_{k+1}\|^2<+\infty$. 
\qed 


Next, we establish the full convergence of the sequence $(x_{k})_{k\in\mathbbm{N}}$ to a solution of problem~\eqref{main-prob}.  The proof is based on the quasi-Fej\'er convergence of the sequence $(x_{k})_{k\in\mathbbm{N}}$ to the set $S_*$, as defined in Definition~\ref{def-fejer}. Note that to establish the convergence of $(x_{k})_{k\in\mathbbm{N}}$ we do not require any local or global Lipschitz continuity assumption of $\nabla f$.

\begin{theorem}[Convergence for the IPG-ELS Method]\label{ptos-de-acum2}
The sequence  $(x_{k})_{k\in\mathbbm{N}}$ generated by the IPG-ELS method converges to a point in
$S_*$. 
\end{theorem}{\it Proof} 
By employing
Proposition~\ref{Prop:Main-proposition}  at $x=x_*\in S_*\subseteq \dom g$, we have
\begin{equation}\label{fjj}
\|x_{k+1}-x_*\|^2\le\|x_{k}-x_*\|^2
+2\left[F(x_{k})-F(x_{k+1})\right]\quad \mbox{for all}\quad k\in \mathbbm{N}.
\end{equation}
We set
$\delta_k:=2\left[F(x_{k})-F(x_{k+1})\right]\ge 0$, and we will prove that $(\delta_k)_{k\in\mathbbm{N}}$ is a summable sequence. In fact,
\begin{align*}
\sum_{k=0}^{+\infty}\delta_k=&
2\sum_{k=0}^ {+\infty}\Big[F(x_{k})-F(x_{k+1})\Big] \le
2\Big[F(x_0)-\lim_{k\to +\infty} F(x_{k+1})\Big]\\ \leq&
2\Big[F(x_0)-F(x_*)\Big]< +\infty.
\end{align*}
 This
together with \eqref{fjj} tells us that the sequence $(x_{k})_{k\in
\mathbbm{N}}$ is quasi-Fej\'er convergent to $S_*$ via
Definition~\ref{def-fejer}. By Lemma \ref{lem:quasi-Fejer}(a), the sequence $(x_{k})_{k\in
\mathbbm{N}}$ is bounded. 
Let $\bar{x}$ be an accumulation point of
$(x_{k})_{k\in\mathbbm{N}}$. Hence, there exists a subsequence $(x_{\ell_{k}})_{k\in
\mathbbm{N}}$ converging to $\bar{x}$.  

Now we proceed by considering the two possible cases:

\noindent {\bf Case 1.} The sequence $\left(\beta_{\ell_{k}}\right)_{k\in
\mathbbm{N}}$ does not converge to $0$, i.e., there exist some
$\beta>0$ and a subsequence of $\left(\beta_{\ell_{k}}\right)_{k\in \mathbbm{N}}$
(without relabeling) such that
\begin{equation}\label{bk-no-0}
\beta_{\ell_{k}}\geq\beta,\quad \forall\, k\in \mathbbm{N}.
\end{equation}By using Proposition~\ref{Prop:Main-proposition} with $x=x_*\in S_*$, we get
\begin{align}\nonumber
\beta_{k}
\left[F(x_{k})-F(x_*)\right]\le&\frac{1}{2}(\|x_{k}-x_*\|^2-\|x_{k+1}-x_*\|^2)
+F(x_{k})-F(x_{k+1}).
\end{align}
 Summing for $k=0, \ldots, m$, the above
inequality implies that
\begin{align*}
\sum_{k=0}^m\beta_{k}
\left[F(x_{k})-F(x_*)\right]\le&\frac{1}{2}(
\|x_0-x_*\|^2-\|x_{m+1}-x_*\|^2)+F(x_0)-F(x_{m+1})\\
\le&\frac{1}{2}\|x_0-x_*\|^2+F(x_0)-F(x_*).
\end{align*}
By taking $m\to+\infty$ and using the fact that $F(x_{k})\ge F(x_*)$ and \eqref{bk-no-0}, we obtain that
$$
\beta\sum_{k=0}^{+\infty}\left[F(x_{\ell_{k}})-F(x_*)\right]\le\sum_{k=0}^{+\infty}\beta_{\ell_{k}}\left[F(x_{\ell_{k}})-F(x_*)\right]\le \sum_{k=0}^{+\infty}\beta_{k}\left[F(x_{k})-F(x_*)\right]<+\infty,
$$
which together with \eqref{bk-no-0} establishes  that
$
\dsty\lim_{k\rightarrow +\infty}\,F(x_{\ell_{k}})=F(x_*).
$ Since $F$ is  lower semicontinuous on $\dom g$, it follows from
the last equality that
\[ F(x_*)\le F(\bar{x})\le
\liminf_{k\rightarrow +\infty}F(x_{\ell_{k}})=\lim_{k\rightarrow +\infty}
F(x_{\ell_{k}})=F(x_*),
\]
 which yields $F(\bar{x})=F(x_*)$
and thus $\bar{x}\in S_*$.

\medskip

\noindent{\bf Case 2.} $\dsty\lim_{k\rightarrow +\infty}\beta_{\ell_{k}}=0$.
Define $\dsty \hat{\beta}_{k}:=\frac{\beta_{k}}{\theta}>0$ and
\begin{equation}\label{2*paso}
\hat{x}_{k+1}:=x_{k}+\hat{\beta}_{k}
(y_{k}-x_{k})=(1-\hat{\beta}_{k})x_{k}+\hat{\beta}_{k}y_{k}, 
\end{equation} where $y_{k}=\tilde{x}_{k}-v_{k}$.
It follows from  the definition of the linesearch
 that
\begin{equation}\label{no-armijo}
f(\hat{x}_{k+1})>f(x_{k})+\hat{\beta}_{k}\la\nabla
f(x_{k}), y_{k}-x_{k}\ra +\frac{\tau}{2}\hat{\beta}_{k}\|x_{k}-\tilde{x}_{k}\|^2+\frac{\gamma_1}{2}\hat{\beta}_{k}\|v_{k}\|^2+\hat{\beta}_{k}\gamma_2\varepsilon_k.
\end{equation} 
It follows from the convexity of $f$, the fact that $\gamma_1>1>\tau$, and the positiveness of the term $\hat{\beta}_{k}\gamma_2\varepsilon_k$ that
\begin{align*}
\la\nabla
f(\hat{x}_{k+1}), \hat{x}_{k+1}-x_{k}\ra \ge f(\hat{x}_{k+1})-f(x_{k})>\hat{\beta}_{k}\la\nabla
f(x_{k}), y_{k}-x_{k}\ra + \hat{\beta}_{k}\frac{\tau}{2}\left (\|x_{k}-\tilde{x}_{k}\|^2+\|v_{k}\|^2\right)\end{align*}
which, together with \eqref{2*paso}, yields 
\begin{align*} \hat{\beta}_{k}\frac{\tau}{2}\left (\|x_{k}-\tilde{x}_{k}\|^2+\|v_{k}\|^2\right)&< \hat{\beta}_{k} \la\nabla
f(\hat{x}_{k+1})-\nabla
f(x_{k}),y_{k}-x_{k}\ra \\ &\le \hat\beta_{k}\| \nabla
f(\hat{x}_{k+1})-\nabla f(x_{k})\|  \| y_{k}-x_{k}\|\\ &= \hat\beta_{k}\| \nabla
f(\hat{x}_{k+1})-\nabla f(x_{k})\|  \| x_{k}-\tilde x_{k}+v_{k}\|.\end{align*}
Now it follows from Lemma \ref{ineq-norm} that
\[
\|x_{k}-\tilde{x}_{k}+v_{k}\|^2 \le 2\left (\|x_{k}-\tilde{x}_{k}\|^2+\|v_{k}\|^2\right).
\]
Hence,
\[
\hat{\beta}_{k}\frac{\tau}{2}\left (\|x_{k}-\tilde{x}_{k}\|^2+\|v_{k}\|^2\right) < \hat{\beta}_{k}\| \nabla f(\hat{x}_{k+1})-\nabla f(x_{k})\| \cdot \sqrt{2}\left (\|x_{k}-\tilde{x}_{k}\|^2+\|v_{k}\|^2\right)^{\frac{1}{2}},
\]
which, due to the positiveness of $\|x_{k}-\tilde{x}_{k}\|^2+\|v_{k}\|^2$, yields
\begin{equation}\label{bb}
\frac{\tau\sqrt{2}}{4}\left(\|x_{k}-\tilde{x}_{k}\|^2+\|v_{k}\|^2\right)^{\frac{1}{2}} \le \| \nabla f(\hat{x}_{k+1})-\nabla f(x_{k})\|.
\end{equation}
Note that $\|\hat{x}_{k+1} - x_{k}\|=\|\hat \beta_k(y_k-x_k)\|=\frac{\hat \beta_k}{\beta_k}\|x_{k+1}-x_k\|=\frac{1}{\theta}\|x_{k+1}-x_k\|,$
which combined with the last statement of Proposition \ref{Prop:Main-proposition} give us that $\|\hat{x}_{\ell_{k}+1} - x_{\ell_{k}}\| \to 0$ as $k \to+\infty$. Since $\nabla f$ is continuous, we have $\| \nabla f(\hat{x}_{\ell_{k}+1})-\nabla f(x_{\ell_{k}})\| \to 0$ as $k \to+\infty$. From \eqref{bb}, it is derived that
\begin{equation}\label{limite1_DI}
\lim_{k \rightarrow+\infty}\, \|x_{\ell_{k}}-\tilde{x}_{\ell_{k}}\|^2 + \|v_{\ell_{k}}\|^2 \le 0,
\end{equation}
therefore, $\lim_{k \rightarrow+\infty}\, \|x_{\ell_{k}}-\tilde{x}_{\ell_{k}}\| = 0$ and $\lim_{k \rightarrow+\infty}\|v_{\ell_{k}}\| = 0$. Additionally, we can use \eqref{ineq-useful-lem-stop} to show that $\lim_{k \rightarrow+\infty}\varepsilon_{\ell_{k}} = 0$. Thus, $\bar{x}$ is also an accumulation point of the sequence $(\tilde{x}_{k})_{k\in \mathbb{N}}$, and $\tilde{x}_{\ell_{k}} \to \bar{x}$ as $k \to+\infty$. Moreover, we have that
\begin{equation}\label{grad-to-0}
\lim_{k \rightarrow+\infty}\| \nabla f(x_{\ell_{k}})-\nabla f(\tilde{x}_{\ell_{k}})\| = 0.
\end{equation}
Now, using \eqref{inex-relative-condition}, we obtain
\[
w_{k} := x_{k} - \tilde{x}_{k} + v_{k} + \nabla f(\tilde{x}_{k}) - \nabla f(x_{k}) \in \nabla f(\tilde{x}_{k}) + \partial_{\varepsilon_k} g(\tilde{x}_{k}) \subseteq \partial_{\varepsilon_k} F(\tilde{x}_{k}).
\]
Furthermore, since $v_{\ell_k}$ converges to $0$ as indicated by \eqref{limite1_DI}, and by applying the triangular inequality, we obtain
\[
\|w_k\| = \|x_k - \tilde{x}_k + v_k + \nabla f(\tilde{x}_k) - \nabla f(x_k) \| \leq \|x_k - \tilde{x}_k\| + \|v_k\| + \|\nabla f(\tilde{x}_k) - \nabla f(x_k)\|,
\]
which implies, via \eqref{limite1_DI} and \eqref{grad-to-0}, that $w_{\ell_k} \in \partial_{\varepsilon_{\ell_k}} F(\tilde{x}_{\ell_k})$ also converges to $0$. Consequently, the convergence of $\tilde{x}_{\ell_k}$ to $\bar{x}$ and $\varepsilon_{\ell_k}$ to $0$, combined with the closedness of the graph of $\partial F$ in Proposition \ref{closed-graph}, gives us that $0 \in \partial F(\bar{x})$. This is equivalent to stating that $\bar{x} \in S_*$.

In all the cases considered above, we have shown that $\bar{x}$, an accumulation point of the sequence $(x_{k})_{k \in \mathbb{N}}$, belongs to $S_*$. Proposition~\ref{lem:quasi-Fejer}(b) implies that $(x_{k})_{k \in \mathbb{N}}$ converges to an optimal solution in $S_*$.
\qed

We continue by showing the convergence rate of the functional values sequence $(F(x_{k}))_{k\in\mathbbm{N}}$.

\begin{theorem}[Convergence Rate of the IPG-ELS Method]\label{theorem:Conv-rate}
Let $(x_{k})_{k\in\mathbbm{N}}$ and $(\beta_{k})_{k\in\mathbbm{N}}$ be generated by the IPG-ELS method. Assume that there exists $\beta > 0$ such that $\beta_{k} \geq \beta$ for all $k \in \mathbbm{N}$. Then, for all $k \in \mathbbm{N}$, we have
\begin{equation}\label{rate-3}
F(x_{k}) - F_* \leq \frac{{\rm dist}(x_0, S_*)^2 + 2\left(F(x_0) - F_*\right)}{2\beta (k+1)}.
\end{equation}
\end{theorem}

{\it Proof} 
For any $\ell\in \mathbbm{N}$ and $x_*\in
S_*$, it follows from  Proposition~\ref{Prop:Main-proposition} with $k=\ell$ and $x=x_*$ that 
\begin{equation*}
F(x_{\ell})-F(x_*)\leq \frac{1}{2\beta_\ell}\left[
\|x_\ell-x_*\|^2-\|x_{\ell+1}-x_*\|^2+2\left(F(x_{\ell})-F(x_{\ell+1})\right)\right]
\end{equation*}
 Summing the above inequality  over $\ell=0,1,\ldots,k$, we have
\begin{align}\label{ineq-sum 0-k}\nonumber
\sum_{\ell=0}^{k}\left(F(x_{\ell})-F(x_*)\right)&\leq\frac{1}{2\beta}\left[
\|x_{0}-x_*\|^2-\|x_{k+1}-x_*\|^2+2\left(F(x_{0})-F(x_{k+1})\right)\right]\\&\leq
\frac{1}{2\beta}\left[
\|x_{0}-x_*\|^2+2\left(F(x_0)-F(x_*)\right)\right].
\end{align}
Since in view of the last statement of Proposition~\ref{Prop:Main-proposition}, we have $F(x_{\ell+1})\leq F(x_{\ell})$ for any  $\ell\in \mathbbm{N}$, it follows from \eqref{ineq-sum 0-k} that 
\begin{equation*}
(k+1)\left(F(x_{k})-F_*\right)\leq
\frac{1}{2\beta}\left(
\|x_0-x_*\|^2+2\left[F(x_{0})-F_*\right]\right).
\end{equation*}
Since  $x_*\in S_*$ is arbitrary, the proof of  \eqref{rate-3} follows.
\qed

\begin{remark}[Complexity of $\eta$-Approximate Solution]
It follows from Theorem~\ref{theorem:Conv-rate} that, given any $\eta>0$, the IPG-ELS method generates an $\eta$-approximate solution $x_{k}$ to problem~\eqref{main-prob}, in the sense that $F(x_{k})-F_*\leq \eta$ in at most $k=\mathcal{O}(1/\eta)$ iterations. We note further that Lemma \ref{lem:Local-Linesearch-Lipschitz} guarantees that the locally Lipschitz assumption for the gradient of $f$ can be used as a sufficient condition for establishing the convergence rate in Theorem~\ref{theorem:Conv-rate}. 
\end{remark}


We end this section by proving the complexity of an $\eta$-approximate stationary solution for problem~\eqref{main-prob} as in Definition~\ref{approxStationarySolution}. For this complexity result, we assume that $\nabla f$ is $L$-Lipschitz continuous on the domain of $g$.  

\begin{theorem}[Complexity of $\eta$-Approximate Stationary Solution]\label{l-rate-456}
Consider $(x_{k})_{k\in\mathbbm{N}}$, $(\tilde x_{k})_{k\in\mathbbm{N}}$, $(v_{k})_{k\in\mathbbm{N}}$, and $(\varepsilon_k)_{k\in\mathbbm{N}}$ generated by the IPG-ELS method and define $w_{k} := v_{k} + x_{k} - \tilde{x}_{k} + \nabla f(\tilde x_{k}) - \nabla f(x_{k})$, for every $k \in \mathbbm{N}$. Then, we have
\begin{equation}\label{inclusion:rk}
w_{k} \in \nabla f(\tilde x_{k}) + \partial_{\varepsilon_k}g(\tilde x_{k}) \subseteq \partial_{\varepsilon_k}F(\tilde x_{k}), \qquad \forall k \in \mathbbm{N}.
\end{equation}
Additionally, if $\alpha \in (0,1-\tau)$, and $\nabla f$ is $L$-Lipschitz continuous on $\dom g$, then, given a tolerance $\eta > 0$, the IPG-ELS method generates an $\eta$-approximate stationary solution $\tilde x_{k}$ to problem~\eqref{main-prob} with residues $(w_{k}, \varepsilon_k)$, in the sense of Definition~\ref{approxStationarySolution}, in at most $k = O\left(1/\eta^2\right)$ iterations.
\end{theorem}

{\it Proof} 
The first inclusion in \eqref{inclusion:rk} follows immediately from \eqref{inexProx-inclusion} and the definition of $w_{k}$, whereas the second inclusion in \eqref{inclusion:rk} follows from the definitions of $F$ and $\partial_\varepsilon F$. 

Now let $x_*$ be the projection of $x_0$ onto $S_*$ and let $d_0 := \|x_0 - x_*\|$. As it was observed in Lemma \ref{lem:Linesearch-Lipschitz}, the Lipschitz continuity of $\nabla f$ implies that there exists $\beta>0$ such that $\beta_{k} \geq \beta$ for all $k$. It follows from Proposition~\ref{Prop:Main-proposition} with $k = \ell \in \mathbbm{N}$ and $x = x_* \in S_*$ that
\begin{equation}
\alpha\beta\|x_{\ell} - \tilde{x}_{\ell}\|^2 \leq \left(\|x_{\ell} - x_*\|^2 - \|x_{\ell+1} - x_*\|^2 + 2\left[F(x_{\ell}) - F(x_{\ell+1})\right]\right) \label{antes-sum-3}
\end{equation}
for all $\ell \in \mathbbm{N}$. Summing the above inequality over $\ell = 0, 1, \ldots, k$, and using that $x_* \in S_*$, we have
\[
\alpha\beta\sum_{\ell=0}^{k}\|x_{\ell} - \tilde{x}_{\ell}\|^2 \leq \|x_0 - x_*\|^2 + 2\left[F(x_0) - F(x_*)\right].
\]
Hence, since $d_0=\|x_0-x_*\|$, we see that there exists ${\ell_{k}}\le k$ such that
\begin{equation}\label{ineq:xktilde-xk}
\|x _{{\ell_{k}}}-\tilde x_{{\ell_{k}}}\|^2\leq \frac{d_0^2+2[F(x_{0})-F(x_*)]}{\alpha\beta (k+1)}.
\end{equation}
On the other hand, since  $\tau, \alpha > 0$,  if follows from \eqref{ineq-useful-lem-stop} that, for every $\ell \in \mathbbm{N}$, 
\begin{align*}
\frac{(\gamma_1-1)}{2}\| v_{\ell}\|^2+\gamma_2 \varepsilon_\ell&\leq \frac{\|x_{\ell}- \tilde x_{\ell}\|^2}{2}+
\langle x_\ell-\tilde x_{\ell},v_{\ell}\rangle\\
&\leq \frac{\|x_{\ell}- \tilde x_{\ell}\|^2}{2}+ \frac{\|x_{\ell}- \tilde x_{\ell}\|^2}{\gamma_1-1}+\frac{(\gamma_1-1)}{4}\|v_{\ell}\|^2,
\end{align*}
where the last inequality is due to Cauchy-Schwarz inequality and the fact that $ab\leq sa^2/2+b^2/(2s)$ for any $a,b\in \mathbbm{R}$ and $s>0$, in particular, with $a=\|v_{\ell}\|$,  $b=\|\tilde x_\ell-x_\ell\|$,  and $s=(\gamma_1-1)/2$. Hence, we have 
\begin{equation}\label{ineq:vkek-yk-xk}
\frac{(\gamma_1-1)}{4}\| v_{\ell}\|^2+\gamma_2 \varepsilon_\ell
\leq \left(\frac{\gamma_1+1}{2(\gamma_1-1)}\right)\|\tilde x_\ell- x_\ell\|^2, \qquad \forall \ell\in \mathbbm{N}.
\end{equation}
Now, from the definition of $w_{\ell}$, \eqref{ineq:vkek-yk-xk}, the Cauchy-Schwarz inequality, the fact that $x_k, \tilde x_k$ are in the domain of $g$,  and the $L$-Lipschitz continuity of $\nabla f$ on the domain of $g$, we have, for every $\ell\in \mathbbm{N}$,
\begin{align}\label{ineq:rk-yk-xk}
\|w_{\ell}\|&\leq \|v_{\ell}\|+\| {x}_{\ell}-\tilde x_\ell\|+\|\nabla f( x_{\ell})-\nabla f( \tilde x_{\ell})\|\nonumber\\ 
& \leq \left[\frac{\sqrt{2(\gamma_1+1)}}{\gamma_1-1}+1+L\right]\|x_{\ell}-\tilde x_{\ell}\|.
\end{align}
 Moreover, it follows from \eqref{ineq:vkek-yk-xk} that 
\begin{equation}\label{ineq:epsk-yk-xk}
\varepsilon_{\ell}
\leq \left(\frac{\gamma_1+1}{2\gamma_2(\gamma_1-1)}\right)\|x_\ell-\tilde x_\ell\|^2.
\end{equation}
Hence, it follows from \eqref{ineq:xktilde-xk},   \eqref{ineq:rk-yk-xk}, and \eqref{ineq:epsk-yk-xk} with $ \ell=\ell_{k}$ and $m_0:=d_0^2+2(F(x_{0})-F(x_*))$ that 
$$
w_{\ell_{k}}\leq \left[\frac{\sqrt{2(\gamma_1+1)}}{\gamma_1-1}+1+L\right]\frac{\sqrt{m_0}}{\sqrt{\alpha\beta (k+1)}}, \qquad \ve_{\ell_{k}} \leq \left(\frac{\gamma_1+1}{2\gamma_2(\gamma_1-1)}\right)\frac{m_0}{\alpha\beta (k+1)}
$$
which in turn implies that
\begin{equation}\label{eq:rlk-epslk=O(k)}
w_{\ell_{k}}=\mathcal{O}(1/\sqrt{k}), \quad \ve_{\ell_{k}}=\mathcal{O}(1/k).
\end{equation}
Thus, the last statement of the theorem follows from \eqref{eq:rlk-epslk=O(k)} and the first inclusion in \eqref{inclusion:rk}.
\qed

\section{Numerical Experiments}\label{NumSec}

In this section, we investigate the numerical behavior of the IPG-ELS method in solving the CUR-like factorization optimization problem \cite{JMLR:v12:mairal11a}. Consider $\mathbbm{E}=\mathbbm{R}^{n \times m}$. Given a matrix $W \in \mathbbm{R}^{m \times n}$, the objective is to find a matrix with sparse rows and columns, $X\in \mathbbm{R}^{n \times m}$, such that $WXW$ approximates $W$. This problem can be formulated as the following splitting optimization problem:
\begin{equation}\label{ee:po09}
\min_{X \in \mathbbm{R}^{n \times m}} \left\{\frac{1}{2}\|W - WXW\|_F^2 + \lambda_{\text{row}}\sum_{i=1}^{n}\|X^{(i)}\|_2 + \lambda_{\text{col}}\sum_{j=1}^{m}\|X_{(j)}\|_2\right\},
\end{equation}
where $\|\cdot\|_F$ denotes the Frobenius norm, and $X^{(i)}$ and $X_{(j)}$ denote the $i$-th row and $j$-th column of $X$, respectively. This problem is a special case of problem~\eqref{main-prob} with
\[
f(X) := \frac{1}{2}\|W - WXW\|_F^2, \quad g(X) := \lambda_{\text{row}}\sum_{i=1}^{n}\|X^{(i)}\|_2 + \lambda_{\text{col}}\sum_{j=1}^{m}\|X_{(j)}\|_2.
\]
In this case, the gradient of $f$ is given by $\nabla f(X) = W^T(WXW-W)W^T$ and has a Lipschitz constant $L = \|W^TW\|_F^2$. 
{   
Note that the proximal operator of \( g \) does not have a closed-form solution; however, \( g \) exhibits a separable structure. By following the approach described in Section~\ref{sep.struct}, we compute a triple \( (\tilde{X}_{k}, 0, \varepsilon_k) \) that satisfies conditions \eqref{inexProx-inclusion}-\eqref{inex-relative-condition} using the Dykstra-like algorithm \cite[Theorem 3.3]{Bauschke:2008c}.
 This algorithm is applied to the proximal subproblem
\begin{equation}\label{eq:675}
\min_{X\in \mathbbm{R}^{n\times m}}  \left\{\frac12\|X-Z\|^2_F+{\lambda_{\rm row}}\sum_{i=1}^{n}\|X^{(i)}\|_2+ {\lambda_{\rm col}}\sum_{j=1}^{m}\|X_{(j)}\|_2\right\},
\end{equation}
where $z:= X_{k}-\nabla f(X_{k}),$ with the initial points $z_0=z, \quad p_0=0, \quad q_0=0$,
 generates the   sequences
\[ 
\left\{ \begin{array}{l}
         y_\ell=\prox_{g_1}(z_\ell+p_\ell)\\
        p_{\ell+1}=z_\ell+p_\ell-y_\ell \end{array} \right. \quad \mbox{and} \quad  \left\{ \begin{array}{l}
         z_{\ell+1}=\prox_{g_2}(y_\ell+q_\ell)\\
        q_{\ell+1}=y_\ell+q_\ell-z_{\ell+1}, \end{array} \right.  
        \]
where $g_1(x)=\lambda_{\rm col}\sum_{j=1}^{m}\|x_{(j)}\|_2$ and $ g_2(x)=\lambda_{\rm row}\sum_{i=1}^{n}\|x^{(i)}\|_2$. 
It follows from the separable property of the proximal operator and Section 6.5.1 of \cite{Parikh-Boyd:2014} that  
$$\prox_{g_1}(x)_{(j)} = \max\left\{1 - \frac{\lambda_{\text{col}}}{\|x_{(j)}\|_2}, 0\right\} x_{(j)},\; \text{for } \quad j=1,\ldots, m$$
and 
$$\prox_{g_2}(x)^{(i)} = \max\left\{1 - \frac{\lambda_{\text{row}}}{\|x^{(i)}\|_2}, 0\right\} x^{(i)},\; \text{for } \quad i=1,\ldots, n.$$
Hence, if the Dykstra-like algorithm is terminated when  
\[
\varepsilon_\ell \leq \frac{(1+\gamma_2)(1-\tau-\alpha)}{2} \|X_k - z_\ell\|^2,
\]  
then $\tilde{X}_k:=z_\ell$ is an inexact proximal solution of \eqref{eq:675} with residual 
\( ( V_k, \varepsilon_k) := ( 0, \varepsilon_\ell) \), where $\varepsilon_\ell$ is as defined in \eqref{eps}.
}

Considering that the IPG-ELS method integrates two effective strategies: (i) permitting the subproblem to be solved inexactly to meet a relative approximation criterion, and (ii) employing an explicit linesearch procedure that computes the proximal operator only once per iteration, our primary goal is to demonstrate that, in certain cases, the proposed method surpasses the proximal gradient method that employs only one of these strategies. { Consequently, we compare the new algorithm with three alternative schemes: an exact version of the IPG-ELS method, denoted by PG-ELS, which corresponds to IPG-ELS with $\gamma_1 = \gamma_2 = 0$, $\theta = 0.5$, $\tau = 1$, and $\epsilon_k\ \leq 10^{-12}$, replacing the inexact criterion in \eqref{inex-relative-condition} with $v_k=0$;  an IPG method with a fixed Stepsize, corresponding to \cite[Algorithm~2]{Millan:2019} with $\alpha_{k} = 1/L$, $\bar w_k=  L(X_{k}-\tilde{X}_{k})- \nabla f(X_{k})$ and $\sigma^2 = 0.9$, where $L$ is the Lipschitz constant of $f$.
 This algorithm is denoted by IPG-FixStep and is defined as $ 
X_{k+1} = \tilde{X}_{k}, $ for all $k \geq 0,$
where $\tilde{X}_{k}$ satisfies
\[
0\in \partial_{\varepsilon_k}g(\tilde{X}_{k}) + L(\tilde{X}_{k} - X_{k}) + \nabla f(X_{k}), \quad   {\varepsilon_k} \leq {0.9L}\|\tilde X_k-\left(X_k-(1/L)\nabla f(x_k)\right)\|^2/2;
\]
and an instance of Tseng's modified forward-backward splitting  method, as described in \cite{Monteiro:2010},  applied to \eqref{ee:po09}. The stepsize is fixed at \( 0.9/L \), and the triple \( (\tilde{X}_{k}, 0, \varepsilon_k) \), associated with the proximal subproblem, is computed to ensure \( \varepsilon_k \leq 10^{-12} \). This algorithm will be referred to as Tseng-MFBS.
}

{
For all methods, the approximate proximal solution $\tilde{X}_{k}$ of \eqref{eq:675}, along with its residual \( (0, \varepsilon_k) \), is computed using the Dykstra-like algorithm as described above.
}
The initialization parameters for the IPG-ELS method were set as $\tau = 0.8$, $\theta = 0.5$, $\gamma_1 = \gamma_2 = 1.1$, and $\alpha = 0.01$. For all tests, we initialized $X_0 = 0 \in \mathbbm{R}^{n \times m}$, and set $\lambda_{\text{row}} = \lambda_{\text{col}} = 0.01$. The IPG-ELS method was executed for $101$ outer iterations to establish a baseline performance metric, $F_* := F(X_{101})$. The other three algorithms were terminated as soon as $F(X_{k}) \leq F_*$ {or after reaching a maximum of 2001 iterations.}. The algorithms were evaluated on six datasets from \cite{Cano:2005,Guyon:2004,Team:2008}: Colon tumor ($62 \times 2000$), heart disease ($303 \times 14$), central nervous system (CNS) ($60 \times 7129$), lung cancer-Michigan ($96 \times 7129$), Secom ($1567 \times 590$), and Cina0 ($132 \times 16033$).

Each matrix $W$ was normalized to have a unit Frobenius norm, with an additional step of centering each column. Subsequently, the resulting matrices were multiplied by a constant, which plays a crucial role in controlling the Lipschitz constant of the function $f$. The experiments were conducted using the Python programming language, which was installed on a machine equipped with a 3.5 GHz Dual-Core Intel Core i5 processor and 16 GB of 2400 MHz DDR4 memory.

In Tables~\ref{table1} and \ref{table2}, we report the Lipschitz constant of the gradient of $f$ (denoted as Lips), the number of outer iterations (O-IT), the number of inner iterations (I-IT), the number of linesearch iterations (LS-IT), and the total running time in seconds (Time). The results indicate that, in terms of CPU times, the IPG-ELS method outperforms the other three methods. This efficiency can be attributed to two main factors: (i) the PG-ELS and Tseng-MFBS methods require  significantly  more inner iterations to solve the proximal subproblem "exactly", and (ii) the IPG-FixStep and Tseng-MFBS methods employ  small stepsizes of $1/L$ and  $0.9/L$, respectively, in the gradient step. 

\begin{table}[H]
\centering
\begin{tabular}{|cccccccc|}
\hline
Problem  & Lips  & Method & $F(X_{k})$ & O-IT& I-IT& LS-IT& Time \\
\hline 
 Colon Tumor & 41.58 & IPG-ELS & 1.1056& 101 & 195& 318 & 25.80 \\ 
 &  & PG-ELS & 1.1056& 104 & 751 & 324& 32.82 \\ 
  &  & IPG-FixStep& 1.1056&1450&1451&  - &  137.78 \\ 
  &  & Tseng-MFBS & 1.1056&1615&4845& - &  193.92 \\ 
 \hline
 Colon Tumor & 665.32 & IPG-ELS & 2.3647& 101 & 101& 823 & 31.02 \\ 
  &  & PG-ELS & 2.3623& 128 & 455 & 1038& 43.40 \\ 
   &  & IPG-FixStep& 2.3647&1103&1104&  - &  94.51 \\ 
   &  & Tseng-MFBS & 2.3645&1227&2454& - &  116.80 \\ 
 \hline
Colon Tumor& 5133.69 & IPG-ELS & 5.8989& 101 & 101& 1134 & 49.30 \\ 
  &  & PG-ELS & 5.8832& 119 & 328 & 1324& 63.73 \\ 
   &  & IPG-FixStep& 5.8981&586&587&  - &  66.26 \\ 
   &  & Tseng-MFBS & 5.8980&652&1304& - &  79.46 \\ 
  \hline 
  Heart Disease& 77.12 & IPG-ELS & 0.1732& 101 & 178& 539 & 1.25 \\ 
  &  & PG-ELS & 0.1732& 119 & 876 & 603& 3.68 \\ 
   &  & IPG-FixStep& 0.1732&487&488&  - &  2.92 \\ 
   &  & Tseng-MFBS & 0.1732&541&1082& - &  4.79 \\  
   \hline
Heart Disease & 1233.99 & IPG-ELS & 0.3129& 101 & 101& 903 & 0.86 \\ 
  &  & PG-ELS & 0.3119& 89 & 368 & 788& 1.71 \\ 
  &  & IPG-FixStep& 0.3520&2001&2001&  - &  11.87 \\ 
  &  & Tseng-MFBS & 0.3577&2001&4002& - &  17.56 \\ 
 \hline
Heart Disease&  9521.56 & IPG-ELS & 0.4995& 101 & 101& 1222 & 1.01 \\ 
 &  & PG-ELS & 0.4992& 227 & 622 & 2690& 3.80 \\ 
   &  & IPG-FixStep& 0.7242&2001&2001&  - &  14.85 \\ 
   &  & Tseng-MFBS & 0.7533&2001&4002& - &  20.74 \\ 
 \hline
  CNS & 41.95 & IPG-ELS & 0.9519& 101 & 182& 341 & 397.09 \\ 
  &  & PG-ELS & 0.9518& 119 & 714 & 396& 531.26 \\ 
   &  & IPG-FixStep& 0.9519&1217&1218&  - &  1762.69 \\ 
   &  & Tseng-MFBS & 0.9519&1354&3091& - &  2263.33 \\ 
 \hline
 CNS&  671.17 & IPG-ELS & 2.1153& 101 & 101& 768 & 554.27 \\ 
 &  & PG-ELS & 2.1119& 148 & 725 & 1110& 885.04 \\ 
   &  & IPG-FixStep& 2.2193&2001&2001&  - &  2752.47 \\ 
   &  & Tseng-MFBS & 2.2439&2001&4002& - &  2828.31 \\ 
 \hline
 CNS&  5178.78 & IPG-ELS & 6.0665& 101 & 101& 1130 & 854.13 \\ 
  &  & PG-ELS & 6.0613& 102 & 343 & 1138& 896.40 \\ 
   &  & IPG-FixStep& 6.0655&775&776&  - &  1240.28 \\ 
   &  & Tseng-MFBS & 6.0653&862&1724& - &  1430.63 \\ 
 \hline
 \end{tabular}
 \caption{Performance of the IPG-ELS, PG-ELS,   IPG-FixStep and Tseng-MFBS methods on $3$ data sets.}
\label{table1}
\end{table}

\begin{table}[H]
\centering
\begin{tabular}{|cccccccc|} \hline
Problem  & Lips  & Method & $F(X_{k})$ & O-IT& I-IT& LS-IT& Time \\
 \hline
 Lung cancer & 52.58 & IPG-ELS & 0.8985& 101 & 179& 443 & 539.37 \\ 
 &  & PG-ELS & 0.8984& 105 & 837 & 457& 646.59 \\ 
   &  & IPG-FixStep& 0.8985&678&679&  - &  1216.98 \\ 
   &  & Tseng-MFBS & 0.8985&755&2265& - &  1487.95 \\ 
 \hline
 Lung cancer & 841.23 & IPG-ELS & 2.7632& 101 & 101& 845 & 601.68 \\ 
  &  & PG-ELS & 2.7623& 131 & 768 & 1085& 875.26 \\ 
   &  & IPG-FixStep& 2.7631&588&589&  - &  850.73 \\ 
   &  & Tseng-MFBS & 2.7631&654&1308& - &  978.97 \\ 
 \hline
  Lung cancer  & 2658.70 & IPG-ELS & 3.6391& 101 & 101& 992 & 711.25 \\ 
  &  & PG-ELS & 3.6378& 161 & 740 & 1574& 1398.46 \\ 
  &  & IPG-FixStep& 3.8819&2001&2001&  - &  3437.99 \\ 
  &  & Tseng-MFBS & 3.9389&2001&4002& - &  3124.96 \\ 
 \hline
Secom & 45.78 & IPG-ELS & 0.6438& 101 & 175& 373 & 99.15 \\ 
  &  & PG-ELS & 0.6438& 99 & 9247 & 360& 822.17 \\ 
   &  & IPG-FixStep& 0.6438&857&858&  - &  304.24 \\ 
   &  & Tseng-MFBS & 0.6438&952&34344& - &  5216.87 \\ 
 \hline
 Secom&  732.51 & IPG-ELS & 0.8587& 101 & 101& 779 & 86.11 \\ 
  &  & PG-ELS & 0.8586& 102 & 4931 & 795& 431.53 \\ 
   &  & IPG-FixStep& 0.8586&1662&1663&  - &  454.59 \\ 
   &  & Tseng-MFBS & 0.8587&1847&3969& - &  637.36 \\ 
 \hline
 Secom&  5652.07 & IPG-ELS & 1.6981& 101 & 101& 1138 & 108.07 \\ 
  &  & PG-ELS & 1.6899& 125 & 1346 & 1388& 220.10 \\ 
   &  & IPG-FixStep& 1.6977&801&802&  - &  218.15 \\ 
   &  & Tseng-MFBS & 1.6977&891&1782& - &  296.50 \\ 
 \hline
 Cina0 & 68.39 & IPG-ELS & 0.7972& 101 & 245& 487 & 3168.99 \\ 
  &  & PG-ELS & 0.7972& 104 & 1490 & 483& 5925.88 \\ 
   &  & IPG-FixStep& 0.7972&609&610&  - &  6451.18 \\ 
   &  & Tseng-MFBS & 0.7972&677&1354& - &  8000.69 \\ 
 \hline
  Cina0& 527.70 & IPG-ELS & 1.2817& 101 & 250& 838 & 4567.13 \\ 
  &  & PG-ELS & 1.2817& 94 & 1789 & 693& 5704.39 \\ 
   &  & IPG-FixStep& 1.2817&1010&1011&  - &  8735.91 \\ 
   &  & Tseng-MFBS & 1.2817&1122&2244& - &  11564.76 \\   
 \hline
 Cina0& 8443.20 & IPG-ELS & 3.6531& 101 & 104& 1126 & 3968.04 \\ 
  &  & PG-ELS & 3.6530& 284 & 3083 & 3061& 15889.08 \\ 
   &  & IPG-FixStep& 3.8493&2001&2001&  - &  17796.07 \\ 
   &  & Tseng-MFBS & 3.8657&2001&4002& - &  20681.39 \\ 
 \hline
 \end{tabular}
 \caption{Performance of the IPG-ELS, PG-ELS,  IPG-FixStep and Tseng-MFBS methods on $3$ data sets.}
\label{table2}
\end{table}

\section{Conclusions}\label{concluding}

In this work, we present an inexact proximal gradient method for solving composite convex optimization problems. This method features a novel explicit linesearch using the relative-type inexact solution of the proximal subproblem. Our approach is primarily designed to solve splitting problems when the objective function is the sum of differentiable and nondifferentiable convex functions, and the analytical computation of the proximal operator is not available. Notably, the convergence of the proposed method is established without assuming Lipschitz continuity of the gradient of the smooth function. This method addresses the need for a balance between computational efficiency and the accuracy of solving the proximal subproblem, a common challenge in practice.

We have confirmed the convergence and iteration complexity of our method, validating its theoretical soundness and practical utility. Numerical experiments demonstrate its applicability and efficiency. Our method maintains convergence rates while efficiently managing relative inexact solutions of the proximal operator. The numerical results indicate that the proposed method competes effectively with both the exact proximal gradient method and the inexact proximal gradient method with a fixed stepsize.

\section*{Declarations}

\noindent{\bf Funding}\\
\textbf{YBC} was partially supported by the NSF Grant DMS-2307328, and by an internal grant from NIU. \textbf{MLNG} was partially supported by CNPq Grants  405349/2021-1 and  304133/2021-3. \textbf{JGM} was partially supported by CNPq Grant 312223/2022-6.

\noindent{\bf Data availability statement}\\
The codes supporting the numerical experiments are freely available in the homepage:\\ https://maxlng.ime.ufg.br/p/17888-publications.

\noindent{\bf Conflicts of interest}\\
The authors declare that they have no conflict of interest.

\begin{acknowledgements} The authors are grateful to the anonymous referees for their valuable comments and suggestions which helped to improve the quality of the
paper.
\end{acknowledgements}


\end{document}